\documentclass[12pt]{amsart}

\usepackage{amsmath}
\usepackage{amssymb}
\usepackage{latexsym}
\usepackage{graphicx} 
\usepackage{mathrsfs}
\usepackage{mathtools}
\usepackage{amsthm}
\usepackage{tikz}
\usetikzlibrary{cd}
\allowdisplaybreaks[4]

\usepackage[colorlinks,citecolor=blue,linkcolor=blue,linktocpage,unicode]{hyperref}

\numberwithin{equation}{section}
\numberwithin{figure}{section}

\newtheorem{thm}{Theorem}[section]
\newtheorem{lem}[thm]{Lemma}
\newtheorem{prop}[thm]{Proposition}
\newtheorem{cor}[thm]{Corollary}

\theoremstyle{definition}
\newtheorem{exam}[thm]{Example}
\newtheorem{defn}[thm]{Definition}
\newtheorem{rem}[thm]{Remark}

\title[Submanifolds with corners in Delzant]{Submanifolds with corners in Delzant polytopes associated to affine subspaces}
\author{Kentaro Yamaguchi}

\date{\today}
\subjclass{Primary 53C40; Secondary 53D20, 14M25}

\keywords{Delzant polytope, manifold with corners, toric K\"{a}hler manifold, torus-equivariantly embedding}

\address{Department of Mathematical Sciences, Tokyo Metropolitan University, 1-1 Minami-Ohsawa, Hachioji, Tokyo, 192-0397, Japan}
\email{yamaguchi-kentaro@ed.tmu.ac.jp}

\begin{document}

\maketitle

\begin{abstract}
In our previous work, we studied the closure of a complex subtorus given from an affine subspace in $\mathfrak{t}^{n} \cong \mathbb{R}^{n}$ in a toric manifold. 
If the closure of the complex subtorus is a smooth complex submanifold in the toric manifold, then we call such submanifold a \textit{torus-equivariantly embedded toric manifold} with respect to the subtorus action determined by the data of the affine subspace.
In this paper, we construct certain submanifolds with corners in a Delzant polytope of the toric manifold from the affine subspaces and show that the conditions for such submanifolds with corners are equivalent to those for torus-equivariantly embedded toric manifolds.
\end{abstract}

\section{Introduction}
\label{sec: introduction}

It is well-known in symplectic toric geometry that there is a bijective correspondence between compact symplectic toric manifolds and Delzant polytopes \cite{MR984900}.
Moreover, for a given Delzant polytope, the corresponding symplectic toric manifold is canonically equipped with the K\"{a}hler structure, known as the Guillemin metric \cite{MR1293656,MR1301331,MR2039163}.
Under the situation, the complements $\check{M}$ of toric divisors in a compact $2n$-dimensional symplectic toric manifold can be identified with an $n$-dimensional complex torus.

In \cite{yamaguchi2023toric},
we study the (Zariski) closure of complex subtori $(\mathbb{C}^{\ast})^{k}$ in a toric manifold $X$ of complex dimension $n$.
More concretely, we consider the properties of the closure (or compactification) $\overline{C(V)}$ of a complex subtori
$C(V) \cong (\mathbb{C}^{\ast})^{k}$ from a $k$-dimensional affine subspace $V$ in $\mathrm{Lie}(T^{n}) = \mathfrak{t}^{n} \cong \mathbb{R}^{n}$.
Note that $\overline{C(V)}$ may not be a complex submanifold in $X$ for some $V$.
If $\overline{C(V)}$ is a smooth submanifold in $X$, then $\overline{C(V)}$ is a torus-equivariantly embedded submanifold with respect to the natural inclusion $i_{V}:T^{k} \to T^{n}$ determined by the data of $V$.
Moreover, the image of the moment map $\overline{\mu}:\overline{C(V)} \to (\mathfrak{t}^{k})^{\ast}$ for the standard Hamiltonian $T^{k}$-action on $\overline{C(V)}$ is equal to the image of the Delzant polytope $\Delta \subset (\mathfrak{t}^{n})^{\ast}$ for $X$ under the pullback $i_{V}^{\ast}: (\mathfrak{t}^{n})^{\ast} \to (\mathfrak{t}^{k})^{\ast}$, i.e., $\overline{\mu}(\overline{C(V)}) = i_{V}^{\ast}(\Delta)$.
We call this submanifold $\overline{C(V)}$ a \textit{torus-equivariantly embedded toric manifold}.

\subsection{Main Results}
In this paper, we further discuss the image of a complex submanifold $\overline{C(V)}$ under the moment map $\mu:X \to (\mathfrak{t}^{n})^{\ast}$ for the $T^{n}$-action on $X$ as a \textit{submanifold} (Definition \ref{def: submanifolds with boundary and corners}) in the Delzant polytope of the ambient toric manifold $X$.
After we see Delzant polytopes $\Delta$ as manifolds with corners (Definition \ref{def: manifold with boundary and corners}, see also \cite{MR3077259} for example), 
we give the conditions of a $k$-dimensional affine subspace $V$ in $\mathbb{R}^{n}$ when $D(V) := \mu(C(V))$ in the interior $\Delta \setminus \partial \Delta$ of the Delzant polytope can be extended to $\Delta$ as a $k$-dimensional submanifold $\overline{D(V)}$ with corners.
In Section \ref{subsec: legendre transformation of affine linear subspaces}, we introduce the notion of \textit{submanifolds with corners} (Definition \ref{def: submanifolds with boundary and corners}).

We show the main result in this paper under the following description:
\begin{equation*}
	\begin{tikzcd}
		\overline{C(V)} \arrow[r,hook] & X \arrow[rr,"\mu"] & & \Delta & \overline{D(V)} \arrow[l, hook] \\
		C(V) \arrow[r, hook] \arrow[rrdd, "\check{\pi}|_{C(V)}"'] \arrow[u, hook] & \check{M} \arrow[rd, "\check{\pi}"'] \arrow[rr, "\mu|_{\check{M}}"] \arrow[u, hook] & &\mathring{\Delta} \arrow[ld, "\Phi"] \arrow[u, hook] & D(V) \arrow[lldd, "\Phi|_{D(V)}"] \arrow[l, hook] \arrow[u, hook] \\
		& & \mathbb{R}^{n} & & \\
		& & V. \arrow[u, hook] & &
	\end{tikzcd}	
\end{equation*}
\begin{thm}[Theorem \ref{thm: lyz correspondence}]
	$\overline{C(V)}$ is a torus-equivariantly embedded toric manifold in a toric manifold $X$ if and only if $\overline{D(V)}$ is a submanifold with corners in a Delzant polytope $\Delta$ of $X$.
\end{thm}

Algebraic geometers have an idea of a correspondence between toric varieties and certain manifolds with corners (see \cite[Section 1.3]{MR0922894} for example).
In this paper, we also deal with this type of the correspondence in a differential-geometric viewpoint after regarding Delzant polytopes as manifolds with corners.
The main theorem might be seen as a generalization of the differential-geometric version of the correspondence to the case of torus-equivariantly embedded toric manifolds.

The construction of $C(V)$ and $D(V)$ is based on Yamamoto's construction \cite{MR3856842} of the Leung--Yau--Zaslow correspondence \cite{MR1894858} in the semi-flat case.
From the Leung--Yau--Zaslow correspondence under the Strominger--Yau--Zaslow mirror symmetry \cite{MR1429831}, it has been proposed that the mirror partners of complex submanifolds in toric manifolds should be Lagrangian submanifolds in the Landau--Ginzbrug mirror of the toric manifolds as conormal bundles of submanifolds in the base manifolds.
Hence, the conormal bundle of $\overline{D(V)}$ should be the mirror partners of $\overline{C(V)}$.

\subsection{Outline}
This paper is organized as follows.  
In Section \ref{sec: preliminary}, we review the settings of \cite{yamaguchi2023toric} and fix the notation. We also show some properties of $\overline{C(V)}$ used in this paper.
In Section \ref{sec: submanifolds in delzant}, 
we introduce the notion of submanifolds with corners 
and give the correspondence between torus-equivariantly embedded toric manifolds and submanifolds with corners in Delzant polytopes.
In Section \ref{subsec: classification of submanifolds in delzant}, we present how to check whether $\overline{D(V)}$ is a submanifold with corners (or equivalently, $\overline{C(V)}$ is a complex submanifold) or not.
In Section \ref{subsec: intersection theory}, we give some remarks regarding the intersection points of torus-equivariantly embedded toric manifolds.

Throughout the paper, we express vectors as column vectors. We identify $(\mathfrak{t}^{n})^{\ast}$ with $\mathbb{R}^{n}$ and use the standard inner product $\langle, \rangle$ on $\mathbb{R}^{n}$.

\section{Reivew of Torus-equivariantly Embedded Toric Manifolds}
\label{sec: preliminary}

In this section, we review the settings in \cite{yamaguchi2023toric}.

\subsection{Toric Manifolds}
\label{subsec: toric}

Due to the Delzant construction \cite{MR984900}, 
compact symplectic toric manifolds are completely classified by \textit{Delzant polytopes}, which are defined as follows:
\begin{defn}
	\label{def: delzant}
	\textit{Delzant polytopes} $\Delta$ are convex polytopes in $(\mathfrak{t}^{n})^{\ast} \cong \mathbb{R}^{n}$ with the following three properties:
	\begin{itemize}
		\item simple; each vertex has $n$ edges,
		\item rational; the direction vectors $v^{\lambda}_{1},\ldots,v^{\lambda}_{n}$ from any vertex $\lambda \in \Lambda$ can be chosen as integral vectors,
		\item smooth; the vectors $v^{\lambda}_{1},\ldots,v^{\lambda}_{n}$ chosen as above form a basis of $\mathbb{Z}^{n}$.
	\end{itemize}
	Here, $\Lambda$ denotes the set of the vertices in $\Delta$.
\end{defn}

Some literature defines Delzant polytopes in terms of the inward pointing normal vectors to each facet instead of the direction vectors from each vertex.
In the rest of paper, we use the following relation between the inward pointing normal vectors and the direction vectors. The proof can be found in \cite[Lemma 3.10]{yamaguchi2023toric} or \cite[Proposition 2.2]{MR2282365}.
\begin{lem}
	\label{lemma: direction vectors vs normal vectors}
	$\Lambda$ denotes the set of the vertices in a Delzant polytope $\Delta$.
	Let $v^{\lambda}_{1},\ldots,v^{\lambda}_{n}$ be the direction vectors from a vertex $\lambda \in \Lambda$ and 
	$u^{\lambda}_{1},\ldots,u^{\lambda}_{n}$ be the inward pointing normal vectors to facets containing the vertex $\lambda \in \Lambda$.
	Then, 
	\begin{equation*}
		\left[ \begin{matrix}
			v^{\lambda}_{1} & \cdots & v^{\lambda}_{n}
		\end{matrix} \right] 
		{^t\! \left[ \begin{matrix}
			u^{\lambda}_{1} & \cdots & u^{\lambda}_{n}
		\end{matrix} \right]}
		= E_{n},
	\end{equation*}
	where $E_{n}$ denotes the identity matrix.
\end{lem}

We define the matrices $Q^{\lambda},D^{\lambda\mu}$ by 
\begin{equation*}
	Q^{\lambda} 
	= \left[ 
		\begin{matrix}
			v^{\lambda}_{1} & \cdots & v^{\lambda}_{n}
		\end{matrix}
	\right]
	(= [Q^{\lambda}_{ij} ]), \;
	D^{\lambda\mu} = (Q^{\lambda})^{-1}Q^{\mu}
	(= [d^{\lambda\mu}_{ij} ]),
\end{equation*}
for any $\lambda,\mu \in \Lambda$.
From Lemma \ref{lemma: direction vectors vs normal vectors}, $\det Q^{\lambda} = \pm 1$. However, we assume that $\det Q^{\lambda} = 1$ by changing the numbering of $v^{\lambda}_{1},\ldots, v^{\lambda}_{n}$.

By using the data of a Delzant polytope, we can see a system of the inhomogeneous coordinate charts on a toric manifold corresponding to the Delzant polytope. 
\begin{lem}
	\label{lemma: coordinate charts}
	From the data of the Delzant polytope $\Delta$ of a toric manifold $X$, 
	we can find an open covering $\{U_{\lambda}\}_{\lambda \in \Lambda}$ of $X$ and a set of maps $\{\varphi_{\lambda}:U_{\lambda} \to \mathbb{C}^{n} \}_{\lambda \in \Lambda}$ such that 
	\begin{equation*}
		\varphi_{\mu}\circ \varphi_{\lambda}^{-1}(z^{\lambda}) 
		= \left(
			\prod_{j=1}^{n}(z^{\lambda}_{j})^{d^{\lambda\mu}_{j1}},
			\ldots,
			\prod_{j=1}^{n}(z^{\lambda}_{j})^{d^{\lambda\mu}_{jn}} 
		\right)
	\end{equation*}
	for any $\lambda,\mu \in \Lambda$ satisfying $U_{\lambda} \cap U_{\mu} \neq \emptyset$.
\end{lem}
The detailed construction is given in \cite[Section 2.1 and Section 3.2]{yamaguchi2023toric}.
Moreover, a K\"{a}hler metric of a toric manifold is determined canonically by the data of the Delzant polytope of the toric manifold, which is called the \textit{Guillemin metric} \cite{MR1293656}.

The complements of toric divisors in a $2n$-dimensional symplectic toric manifold can be identified with a complex $n$-dimensional torus. 
Using the data of a Delzant polytope, the identification is written explicitly (see \cite[Section 3.3]{yamaguchi2023toric} for detail in terms of the notation used here).

\subsection{Torus-equivariantly Embedded Toric Manifolds}
\label{subsec: toric submanifolds}
Let $V$ be an affine subspace in $\mathbb{R}^{n}$ and $p_{1},\ldots, p_{k} \in \mathbb{Z}^{n}$ be primitive vectors which form a basis of the linear part of $V$.
Thus, $V$ can be written as $V = \mathbb{R}p_{1} + \cdots + \mathbb{R}p_{k} + a$ for $a \in \mathbb{R}^{n}$. Note that $V$ should have rational slope.

As Yamamoto considered in \cite{MR3856842}, a $k$-dimensional complex subtorus $C(V)$ in $(\mathbb{C}^{\ast})^{n}$ can be constructed from the data of $V$.
\begin{prop}[{\cite[Proposition 4.2]{yamaguchi2023toric}}]
	\label{prop: complex submanifolds in complex torus2}
	Given a $k$-dimensional affine subspace $V$ in $\mathbb{R}^{n}$,
	there exists a primitive basis $q_{k+1}, \ldots, q_{n} \in \mathbb{Z}^{n}$ of the orthogonal subspace to $V = \mathbb{R}p_{1} + \cdots + \mathbb{R}p_{k} + a$ in $\mathbb{R}^{n}$ such that 
	\begin{equation*}
		C(V) = \left\{ (e^{w_{1}}, \ldots, e^{w_{n}}) \in (\mathbb{C}^{\ast})^{n} \;\middle|\; ~^t[q_{k+1} \cdots q_{n}] \left( \left[ \begin{array}{c} w_{1} \\ \vdots \\ w_{n} \end{array} \right] -a \right) = 0 \right\},
	\end{equation*}
	where $w_{1} = x_{1} + \sqrt{-1}y_{1}, \ldots, w_{n} = x_{n} + \sqrt{-1}y_{n}$.
\end{prop}

Recall that the complements of toric divisors can be seen as a complex torus. 
By taking the compactification (or taking the closure) $\overline{C(V)}$ of $C(V) \cong (\mathbb{C}^{\ast})^{k}$ in a toric manifold, we may find a submanifold of the toric manifold explicitly as below.
\begin{defn}
	\label{def: compactification}
	\textit{The compactification $\overline{C(V)} \subset X$ of $C(V)$} is defined as $\overline{C(V)} = \bigcup_{\lambda \in \Lambda}\varphi_{\lambda}^{-1}(\overline{C_{\lambda}(V)})$, where
	\begin{equation*}
		\overline{C_{\lambda}(V)} 
		= \{ z^{\lambda} \in \varphi_{\lambda}(U_{\lambda}) \mid f^{\lambda}_{j}(z^{\lambda}) = 0, j= k+1,\ldots,n \}.
	\end{equation*}
	Here, the map $f^{\lambda}_{j}$ is defined by
	\begin{equation}
		\label{eq: defining eqautions for toric submanifolds}
		f^{\lambda}_{j}(z^{\lambda}) 
		= \prod_{i \in \mathcal{I}^{+}_{\lambda,j}} (z^{\lambda}_{i})^{\langle u^{\lambda}_{i}, q_{j} \rangle}
		- e^{\langle a, q_{j} \rangle}\prod_{i \in \mathcal{I}^{-}_{\lambda,j}} (z^{\lambda}_{i})^{-\langle u^{\lambda}_{i}, q_{j} \rangle} 
	\end{equation}
	for each $j = k+1, \ldots, n$, where 
	\begin{align}
		\label{eq: plus}
		\mathcal{I}^{+}_{\lambda,j} =& \left\{ i \in \{1,2,\ldots, n \} \;\middle|\; \langle u^{\lambda}_{i}, q_{j} \rangle \geq 0 \right\}, \\ 
		\label{eq: minus}
		\mathcal{I}^{-}_{\lambda,j} =& \left\{ i \in \{1,2,\ldots, n \} \;\middle|\; \langle u^{\lambda}_{i}, q_{j} \rangle \leq 0 \right\}, \\
		\mathcal{I}^{0}_{\lambda,j} =& \left\{ i \in \{1,2,\ldots, n \} \;\middle|\; \langle u^{\lambda}_{i}, q_{j} \rangle = 0 \right\}.
	\end{align}
	Here, if $\mathcal{I}^{+}_{\lambda,j}\setminus \mathcal{I}^{0}_{\lambda,j} = \emptyset$, then $\prod_{i \in \mathcal{I}^{+}_{\lambda,j}} (z^{\lambda}_{i})^{\langle u^{\lambda}_{i}, q_{j} \rangle} = 1$. Similarly, if $\mathcal{I}^{-}_{\lambda,j} \setminus \mathcal{I}^{0}_{\lambda,j}= \emptyset$, then $\prod_{i \in \mathcal{I}^{-}_{\lambda,j}} (z^{\lambda}_{i})^{-\langle u^{\lambda}_{i}, q_{j} \rangle} = 1$.
\end{defn}

In fact, $\overline{C_{\lambda}(V)}$ is a zero locus of $f^{\lambda}_{k+1},\ldots, f^{\lambda}_{n}$.
By the implicit function theorem, $\overline{C(V)}$ is a complex submanifold in $X$ if the rank of the Jacobian matrix of $f^{\lambda}:= [f^{\lambda}_{k+1}\; f^{\lambda}_{k+2} \;\cdots \; f^{\lambda}_{n}]$ is equal to $n-k$ for any points $p \in \overline{C(V)}$ and any $\lambda \in \Lambda$.

\begin{rem}
	\label{remark: toric submanifold}
	It is obvious that $C(V)$ is always a complex submanifold in $(\mathbb{C}^{\ast})^{n}$ for any affine subspace $V$ in $\mathbb{R}^{n}$, 
	but $\overline{C(V)}$ may not be a complex submanifold in a toric manifold $X$.
	In \cite[Section 5]{yamaguchi2023toric}, we gave examples of affine subspaces $V$ and toric manifolds $X$ for $\overline{C(V)}$ not to be complex submanifolds in $X$.
\end{rem}

In \cite[Lemma 4.14]{yamaguchi2023toric}, it was shown that the vectors $p_{1},\ldots, p_{k}$, $q_{k+1},\ldots, q_{n}$, $v^{\lambda}_{1}, \ldots, v^{\lambda}_{n}$, and $u^{\lambda}_{1},\ldots,u^{\lambda}_{n}$ have the following relation:
\begin{lem}
	For any $l=1,\ldots,k$ and any $j=k+1,\ldots, n$, 
	\begin{equation*}
		\sum_{i=1}^{n}\langle p_{l},v^{\lambda}_{i} \rangle \langle u^{\lambda}_{i},q_{j} \rangle = 0
	\end{equation*}
	holds.
\end{lem}

We review the definition of the $T^{k}$-action on $\overline{C(V)}$ given in \cite[Section 4.3]{yamaguchi2023toric}.
Recall that we defined the map $i_{V}:T^{k} \to T^{n}$ by 
\begin{equation*}
	i_{V}(t_{1},\ldots, t_{k})
	=
	\left(
		\prod_{l=1}^{k} t_{l}^{\langle p_{l},e_{1} \rangle}, \ldots, \prod_{l=1}^{k} t_{l}^{\langle p_{l},e_{n} \rangle}
	\right).
\end{equation*}
Here, $e_{1},\ldots,e_{n}$ denote the standard basis of $\mathbb{R}^{n}$.
We can describe the $T^{k}$-action on $\varphi_{\lambda}(U_{\lambda})$ by 
\begin{align*}
	T^{k} \times \varphi_{\lambda}(U_{\lambda}) &\to \varphi_{\lambda}(U_{\lambda}) \\
	(t = (t_{1},\ldots, t_{k}), z^{\lambda}) &\mapsto i_{V}(t)\cdot z^{\lambda},
\end{align*}
where $i_{V}(t)\cdot z^{\lambda}$ means that 
\begin{equation*}
	i_{V}(t)\cdot z^{\lambda} 
	= 
	\left(
		\prod_{l=1}^{k} t_{l}^{\langle p_{l},v^{\lambda}_{1} \rangle}z^{\lambda}_{1}, \ldots, \prod_{l=1}^{k} t_{l}^{\langle p_{l},v^{\lambda}_{n} \rangle}z^{\lambda}_{n}
	\right).
\end{equation*}
Since we have a natural inclusion $i_{\lambda}:\overline{C_{\lambda}(V)} \to \varphi_{\lambda}(U_{\lambda})$, we can define the $T^{k}$-action on $\overline{C_{\lambda}(V)}$ by 
\begin{align*}
	T^{k} \times \overline{C_{\lambda}(V)} &\to \overline{C_{\lambda}(V)} \\
	(t, z^{\lambda}) &\mapsto i_{V}(t)\cdot z^{\lambda},
\end{align*}
which makes the following diagram commutative:
\begin{equation*}
	\begin{tikzcd}
		T^{n} & \times & \varphi_{\lambda}(U_{\lambda}) \arrow[r] & \varphi_{\lambda}(U_{\lambda}) \\
		T^{k} \arrow[u, "i_{V}"'] & \times & \overline{C_{\lambda}(V)} \arrow[r] \arrow[u, "i_{\lambda}"'] & \overline{C_{\lambda}(V)} \arrow[u, "i_{\lambda}"'].
	\end{tikzcd}
\end{equation*}

\subsection{Properties of the Rank of the Jacobian Matrices}
\label{subsec: properties on the rank of the jacobian matrix}

In this section, we show properties of the rank of the Jacobian matrix $Df^{\lambda}$ of $f^{\lambda}= [f^{\lambda}_{k+1}\; \cdots \; f^{\lambda}_{n}]$.

We first show that the rank of 
the Jacobian matrix of $f^{\lambda}= [f^{\lambda}_{k+1}\; \cdots \; f^{\lambda}_{n}]$ defined in Equation (\ref{eq: defining eqautions for toric submanifolds}) is invariant under the $T^{k}$-action on $\overline{C(V)}$ for each $\lambda \in \Lambda$.

\begin{lem}
	\label{lemma: subtorus action on defining equations}
	For $t = (t_{1},\ldots,t_{k}) \in T^{k}$, we have 
	\begin{equation}
		\label{eq: subtorus action on defining equations}
		f^{\lambda}_{i}(i_{V}(t) \cdot z^{\lambda}) = \mathcal{T}^{\lambda}_{i}(t) f^{\lambda}_{i}(z^{\lambda})
	\end{equation}
	for $i = k+1,\ldots, n$, where the functions $f^{\lambda}_{k+1},\ldots,f^{\lambda}_{n}$ are defined in Equation (\ref{eq: defining eqautions for toric submanifolds}) and $\mathcal{T}^{\lambda}_{k+1}(t), \ldots, \mathcal{T}^{\lambda}_{n}(t)$ are defined by 
	\begin{equation*}
		\mathcal{T}^{\lambda}_{i}(t)
		= \prod_{j\in \mathcal{I}^{-}_{\lambda,i}} \prod_{l=1}^{k} t_{l}^{-\langle p_{l},v^{\lambda}_{j} \rangle \langle u^{\lambda}_{j}, q_{i} \rangle}
	\end{equation*}
	for $i = k+1,\ldots, n$.

\begin{proof}
	Since we define the $T^{k}$-action on $\varphi_{\lambda}(U_{\lambda})$ in Section \ref{subsec: toric submanifolds}, we write $f^{\lambda}_{i}(i_{V}(t) \cdot z^{\lambda})$ for $i=k+1,\ldots, n$ as 
	\begin{align*}
		f^{\lambda}_{i}(i_{V}(t) \cdot z^{\lambda})
		=&
		\prod_{j \in \mathcal{I}^{+}_{\lambda,i}} 
		\left(
			\prod_{l=1}^{k} t_{l}^{\langle p_{l}, v^{\lambda}_{j} \rangle} z^{\lambda}_{j} 
		\right)^{\langle u^{\lambda}_{j},q_{i} \rangle} \\
		&- 
		e^{\langle a,q_{i} \rangle}
		\prod_{j \in \mathcal{I}^{-}_{\lambda,i}} 
		\left(
			\prod_{l=1}^{k} t_{l}^{\langle p_{l}, v^{\lambda}_{j} \rangle} z^{\lambda}_{j} 
		\right)^{-\langle u^{\lambda}_{j},q_{i} \rangle} \\
		=&
		\prod_{j \in \mathcal{I}^{+}_{\lambda,i}} 
		\prod_{l=1}^{k}
		t_{l}^{\langle p_{l},v^{\lambda}_{j} \rangle \langle u^{\lambda}_{j}, q_{i} \rangle} 
		(z^{\lambda}_{j})^{\langle u^{\lambda}_{j}, q_{i} \rangle} \\
		&- 
		e^{\langle a,q_{i} \rangle}
		\prod_{j \in \mathcal{I}^{-}_{\lambda,i}} 
		\prod_{l=1}^{k}
		t_{l}^{-\langle p_{l},v^{\lambda}_{j} \rangle \langle u^{\lambda}_{j}, q_{i} \rangle} 
		(z^{\lambda}_{j})^{-\langle u^{\lambda}_{j}, q_{i} \rangle}.
	\end{align*}
	Since $\mathcal{I}^{+}_{\lambda,i} \cup \mathcal{I}^{-}_{\lambda,i} = \{1,\ldots, n \}$, we calculate
	\begin{align*}
		\frac{\displaystyle\prod_{j \in \mathcal{I}^{+}_{\lambda,i}} 
		\prod_{l=1}^{k}
		t_{l}^{\langle p_{l},v^{\lambda}_{j} \rangle \langle u^{\lambda}_{j}, q_{i} \rangle}}{\displaystyle\prod_{j \in \mathcal{I}^{-}_{\lambda,i}} \prod_{l=1}^{k}
		t_{l}^{-\langle p_{l},v^{\lambda}_{j} \rangle \langle u^{\lambda}_{j}, q_{i} \rangle}}
		=&
		\prod_{l=1}^{k} t_{l}^{\sum_{j \in \mathcal{I}^{+}_{\lambda,i}} \langle p_{l},v^{\lambda}_{j} \rangle \langle u^{\lambda}_{j}, q_{i} \rangle + \sum_{j \in \mathcal{I}^{-}_{\lambda,i}}\langle p_{l},v^{\lambda}_{j} \rangle \langle u^{\lambda}_{j}, q_{i} \rangle} \\
		=&
		\prod_{l=1}^{k}
		t_{l}^{\sum_{j=1}^{n} \langle p_{l},v^{\lambda}_{j} \rangle \langle u^{\lambda}_{j}, q_{i} \rangle} \\
		=&
		1.
	\end{align*} 
	Thus, we have 
	\begin{align*}
		f^{\lambda}_{i}(i_{V}(t) \cdot z^{\lambda})
		=&
		\prod_{j \in \mathcal{I}^{+}_{\lambda,i}} 
		\prod_{l=1}^{k}
		t_{l}^{\langle p_{l},v^{\lambda}_{j} \rangle \langle u^{\lambda}_{j}, q_{i} \rangle} 
		(z^{\lambda}_{j})^{\langle u^{\lambda}_{j}, q_{i} \rangle} \\
		&- 
		e^{\langle a,q_{i} \rangle}
		\prod_{j \in \mathcal{I}^{-}_{\lambda,i}} 
		\prod_{l=1}^{k}
		t_{l}^{-\langle p_{l},v^{\lambda}_{j} \rangle \langle u^{\lambda}_{j}, q_{i} \rangle} 
		(z^{\lambda}_{j})^{-\langle u^{\lambda}_{j}, q_{i} \rangle} \\
		=&
		\mathcal{T}^{\lambda}_{i}(t)
		\left(
			\prod_{j \in \mathcal{I}^{+}_{\lambda,i}}
			(z^{\lambda}_{j})^{\langle u^{\lambda}_{j}, q_{i} \rangle}
			-
			e^{\langle a,q_{i} \rangle}
			\prod_{j \in \mathcal{I}^{-}_{\lambda,i}} 
			(z^{\lambda}_{j})^{-\langle u^{\lambda}_{j}, q_{i} \rangle}
		\right) \\
		=&
		\mathcal{T}^{\lambda}_{i}(t)
		f^{\lambda}_{i}(z^{\lambda}),
	\end{align*}
	which is Equation (\ref{eq: subtorus action on defining equations}).
\end{proof}
\end{lem}

\begin{prop}
	\label{prop: subtorus action on Jacobian matrix}
	In the above situation, we have 
	\begin{equation}
		\label{eq: subtorus action on Jacobian matrix}
		Df^{\lambda}(i_{V}(t) \cdot z^{\lambda})
		= \mathcal{T}^{\lambda}(t) Df^{\lambda}(z^{\lambda}) \mathcal{T}^{\lambda}_{V}(t),
	\end{equation}
	where the matrices $\mathcal{T}^{\lambda}(t)$ and $\mathcal{T}^{\lambda}_{V}$ are defined by 
	\begin{align*}
		\mathcal{T}^{\lambda}(t)
		= 
		\left[
			\begin{matrix}
				\mathcal{T}^{\lambda}_{k+1}(t) & & \\
				& \ddots & \\
				& & \mathcal{T}^{\lambda}_{n}(t)
			\end{matrix}
		\right], \;
		\mathcal{T}^{\lambda}_{V} (t)
		= 
		\left[
			\begin{matrix}
				\displaystyle
				\prod_{l=1}^{k}t_{l}^{-\langle p_{l},v^{\lambda}_{1} \rangle} & & \\
				& \ddots & \\
				& & \displaystyle \prod_{l=1}^{k}t_{l}^{-\langle p_{l},v^{\lambda}_{n} \rangle}
			\end{matrix}
		\right].
	\end{align*}
	In particular, the rank of $Df^{\lambda}(i_{V}(t)\cdot z^{\lambda})$ is equal to the rank of $Df^{\lambda}(z^{\lambda})$.

\begin{proof}
	From Equation (\ref{eq: subtorus action on defining equations}), we have 
	\begin{align*}
		\frac{\partial f^{\lambda}_{i}}{\partial z^{\lambda}_{j}}(i_{V}(t)\cdot z^{\lambda})
		=&
		\lim_{h^{\prime} \to 0}\frac{f^{\lambda}_{i}(i_{V}(t)\cdot z^{\lambda} - h^{\prime}e_{j}) - f^{\lambda}_{i}(i_{V}(t)\cdot z^{\lambda})}{h^{\prime}} \\
		=&
		\lim_{h \to 0}\frac{f^{\lambda}_{i}(i_{V}(t)\cdot (z^{\lambda} - he_{j})) - f^{\lambda}_{i}(i_{V}(t)\cdot z^{\lambda})}{\prod_{l=1}^{k}t_{l}^{\langle p_{l},v^{\lambda}_{j}\rangle}h} \\
		=&
		\mathcal{T}^{\lambda}_{i}(t)
		\prod_{l=1}^{k}t_{l}^{-\langle p_{l},v^{\lambda}_{j}\rangle}
		\lim_{h \to 0}\frac{f^{\lambda}_{i}(z^{\lambda} - he_{j}) - f^{\lambda}_{i}(z^{\lambda})}{h} \\
		=&
		\mathcal{T}^{\lambda}_{i}(t)
		\prod_{l=1}^{k}t_{l}^{-\langle p_{l},v^{\lambda}_{j}\rangle}
		\frac{\partial f^{\lambda}_{i}}{\partial z^{\lambda}_{j}}(z^{\lambda}).
	\end{align*}
	Thus, we can express the matrix $Df^{\lambda}(i_{V}(t)\cdot z^{\lambda})$ as 
	\begin{align*}
		&Df^{\lambda}(i_{V}(t)\cdot z^{\lambda}) \\
		=&
		\left[
			\begin{matrix}\displaystyle
				\mathcal{T}^{\lambda}_{k+1}(t)
				\prod_{l=1}^{k}t_{l}^{-\langle p_{l},v^{\lambda}_{1}\rangle}
				\frac{\partial f^{\lambda}_{k+1}}{\partial z^{\lambda}_{1}}(z^{\lambda}) 
				& \cdots & \displaystyle
				\mathcal{T}^{\lambda}_{k+1}(t)
				\prod_{l=1}^{k}t_{l}^{-\langle p_{l},v^{\lambda}_{n}\rangle}
				\frac{\partial f^{\lambda}_{k+1}}{\partial z^{\lambda}_{n}}(z^{\lambda}) \\
				\vdots & \ddots & \vdots \\ \displaystyle
				\mathcal{T}^{\lambda}_{n}(t)
				\prod_{l=1}^{k}t_{l}^{-\langle p_{l},v^{\lambda}_{1}\rangle}
				\frac{\partial f^{\lambda}_{n}}{\partial z^{\lambda}_{1}}(z^{\lambda}) 
				& \cdots & \displaystyle
				\mathcal{T}^{\lambda}_{n}(t)
				\prod_{l=1}^{k}t_{l}^{-\langle p_{l},v^{\lambda}_{n}\rangle}
				\frac{\partial f^{\lambda}_{n}}{\partial z^{\lambda}_{n}}(z^{\lambda})
			\end{matrix}
		\right] \\
		=&
		\mathcal{T}^{\lambda}(t)
		Df^{\lambda}(z^{\lambda})
		\mathcal{T}^{\lambda}_{V}(t).
	\end{align*}
	Therefore, we obtain Equation (\ref{eq: subtorus action on Jacobian matrix}).

	Since the matrices $\mathcal{T}^{\lambda}(t)$ and $\mathcal{T}^{\lambda}_{V}(t)$ are both invertible for any $t \in T^{k}$, we conclude that the rank of the matrix $Df^{\lambda}(z^{\lambda})$ is equal to the rank of the matrix $Df^{\lambda}(i_{V}(t)\cdot z^{\lambda})$ for any $\lambda \in \Lambda$ and $t \in T^{k}$.
\end{proof}
\end{prop}

We show that the rank of our Jacobian matrix $Df^{\lambda}$ at a point $z^{\lambda} \in \overline{C_{\lambda}(V)}$ is determined by the absolute value of $z^{\lambda}$.

\begin{lem}
	\label{lemma: subtorus action to absolute value}
	For any point $z^{\lambda} = (z^{\lambda}_{1},\ldots, z^{\lambda}_{n}) \in \overline{C_{\lambda}(V)}$, there exists $t \in T^{k}$ such that 
	\begin{equation*}
		i_{V}(t)\cdot z^{\lambda} = (\vert z^{\lambda}_{1}\vert, \ldots, \vert z^{\lambda}_{n} \vert ).
	\end{equation*}
\begin{proof}
	If $z^{\lambda}_{i} = 0$, then $\vert z^{\lambda}_{i} \vert = z^{\lambda}_{i} = 0$. We assume that $z^{\lambda}_{i} \neq 0$.

	We use the expression given in \cite[Lemma 4.1]{yamaguchi2023toric}, i.e., if $(z_{1},\ldots, z_{n})\in C(V)$, then for any $i=1,\ldots,n$,
	\begin{equation*}
		z_{i} = \exp \left(\sum_{l=1}^{k}\langle p_{l},e_{i} \rangle u_{l} + \sqrt{-1}\langle a,e_{i} \rangle + \sqrt{-1}\sum_{l=1}^{k} \langle p_{l},e_{i} \rangle v_{l} \right),
	\end{equation*}
	where $(e^{u_{1}+\sqrt{-1}v_{1}},\ldots,e^{u_{k}+\sqrt{-1}v_{k}}) \in (\mathbb{C}^{\ast})^{k}$.
	In \cite[Definition 3.16]{yamaguchi2023toric}, we defined the map $\hat{\phi}_{\lambda}:(\mathbb{C}^{\ast})^{n} \to \varphi_{\lambda}(U_{\lambda})$ by 
	\begin{equation*}
		\hat{\phi}_{\lambda}(z_{1},\ldots,z_{n})
		=
		\left(
			\prod_{j=1}^{n}(z_{j})^{\langle e_{j} v^{\lambda}_{1} \rangle}, \ldots, \prod_{j=1}^{n}(z_{j})^{\langle e_{j} v^{\lambda}_{n} \rangle}
		\right).
	\end{equation*}
	We obtain that the $i$-th element $z^{\lambda}_{i}$ of $\hat{\phi}_{\lambda}(z_{1},\ldots,z_{n})$ is 
	\begin{align*}
			z^{\lambda}_{i} = \exp \left( \sum_{l=1}^{k}\langle p_{l},v^{\lambda}_{i}\rangle u_{l} + \langle a,v^{\lambda}_{i} \rangle + \sqrt{-1}\sum_{l=1}^{k}\langle p_{l},v^{\lambda}_{i}\rangle v_{l} \right).
	\end{align*}
	We define $t \in T^{k}$ by $t_{l} = e^{\sqrt{-1}\theta_{l}}$ ($\theta_{1},\ldots, \theta_{k} \in \mathbb{R}{/ 2\pi \mathbb{Z}}$).
	By the $T^{k}$-action, we obtain 
	\begin{align*}
		&\prod_{l=1}^{k}t_{l}^{\langle p_{l},v^{\lambda}_{i} \rangle}z^{\lambda}_{i} \\
		&=
		\exp \left(\sum_{l=1}^{k} \langle p_{l},v^{\lambda}_{i} \rangle u_{l} + \langle a,v^{\lambda}_{i} \rangle \right) \exp \left( \sqrt{-1}\sum_{l=1}^{k} \langle p_{l},v^{\lambda}_{i} \rangle (v_{l} + \theta_{l}) \right).
	\end{align*}
	For any $v_{1},\ldots, v_{k}$, we can define $\theta_{1} = -v_{1},\ldots, \theta_{k} = -v_{k}$. Therefore, we obtain that for such $(t_{1},\ldots,t_{k}) \in T^{k}$,
	\begin{align*}
		&\prod_{l=1}^{k}t_{l}^{\langle p_{l},v^{\lambda}_{i} \rangle}z^{\lambda}_{i} \\
		&=
		\prod_{l=1}^{k}(e^{-\sqrt{-1}v_{l}})^{\langle p_{l},v^{\lambda}_{i} \rangle}\exp \left( \sum_{l=1}^{k}\langle p_{l},v^{\lambda}_{i}\rangle u_{l} + \langle a,v^{\lambda}_{i} \rangle + \sqrt{-1}\sum_{l=1}^{k}\langle p_{l},v^{\lambda}_{i}\rangle v_{l} \right) \\
		&= 
		\exp \left( \sum_{l=1}^{k}\langle p_{l},v^{\lambda}_{i}\rangle u_{l} + \langle a,v^{\lambda}_{i} \rangle \right) \\
		&= 
		\vert z^{\lambda}_{i} \vert,
	\end{align*}
	which is the desired result.
\end{proof}
\end{lem}

By Proposition \ref{prop: subtorus action on Jacobian matrix} and Lemma \ref{lemma: subtorus action to absolute value}, we obtain the following:

\begin{cor}
	\label{cor: subtorus action on Jacobian matrix2}
	For any point $z^{\lambda} \in \overline{C_{\lambda}(V)}$, the rank of the matrix $Df^{\lambda}(z^{\lambda})$ is equal to the rank of the matrix $Df^{\lambda}(\vert z^{\lambda} \vert)$, where $\vert z^{\lambda} \vert$ means that $\vert z^{\lambda} \vert = (\vert z^{\lambda}_{1} \vert, \ldots, \vert z^{\lambda}_{n} \vert)$.
\end{cor}

We end this section by showing that if $\dim V = n-1$, then the rank of the Jacobian matrix $Df^{\lambda}$ is determined only by the linear part of the affine subspace $V$.
\begin{prop}
\label{prop: affine not affect when codim1}
Let $V = \mathbb{R}p_{1}+\cdots +\mathbb{R}p_{k} + a \; (a \in \mathbb{R}^{n})$ be an affine subspace in $\mathbb{R}^{n}$.
If $n-k = 1$, then the rank of the Jacobian matrix of $f^{\lambda}$ is independent of the choice of $a \in \mathbb{R}^{n}$ for any $\lambda \in \Lambda$.
\begin{proof}
	We define $C^{\lambda}_{j}$ by 
	\begin{equation*}
		C^{\lambda}_{j} =
		\begin{cases*}
			0 & if $j \in \mathcal{I}^{+}_{\lambda}$, \\
			1 & if $j \in \mathcal{I}^{-}_{\lambda} \setminus \mathcal{I}^{0}_{\lambda}$,
		\end{cases*}
	\end{equation*}
	for any $j = 1,\ldots, n$, and $\lambda \in \Lambda$.

	For any $\lambda \in \Lambda$, the Jacobian matrix $Df^{\lambda}$ is written as 
	\begin{align*}
		Df^{\lambda}
		=&
		\left[
			\begin{matrix}
				\displaystyle
				\frac{\partial f^{\lambda}}{\partial z^{\lambda}_{1}} 
				& \cdots & 
				\displaystyle
				\frac{\partial f^{\lambda}}{\partial z^{\lambda}_{n}}
			\end{matrix}
		\right] \\
		=&
		\left[
			\begin{matrix}
				\displaystyle
				\langle u^{\lambda}_{1},q \rangle \prod_{l \in \mathcal{I}^{\pm}_{\lambda}} (z^{\lambda}_{l})^{\pm \langle u_{l},q \rangle - \delta_{l1}} 
				& \cdots &
				\displaystyle
				\langle u^{\lambda}_{n},q \rangle \prod_{l \in \mathcal{I}^{\pm}_{\lambda}} (z^{\lambda}_{l})^{\pm \langle u_{l},q \rangle - \delta_{ln}} 
			\end{matrix}
		\right] \tilde{C}^{\lambda}_{a},
	\end{align*}
	where $\tilde{C}^{\lambda}_{a}$ is given by
	\begin{equation*}
		\tilde{C}^{\lambda}_{a}
		=
		\left[
			\begin{matrix}
				e^{\langle a,q \rangle C^{\lambda}_{1}} & &\\
				& \ddots & \\
				& & e^{\langle a,q \rangle C^{\lambda}_{n}}
			\end{matrix}
		\right].
	\end{equation*}
	Since the matrix $\tilde{C}^{\lambda}_{a}$ is invertible, we see that the rank of $Df^{\lambda}$ is independent of the choice of $a \in \mathbb{R}^{n}$.
\end{proof}
\end{prop}

At this moment, the author does not know whether the rank of the Jacobian matrix $Df^{\lambda}$ is independent of the choice of $a \in \mathbb{R}^{n}$ for any $\lambda \in \Lambda$ when $n-k > 1$.

\section{Submanifolds in Delzant Polytopes}
\label{sec: submanifolds in delzant}

In this section, we consider the image of the moment map of a torus-equivariantly embedded toric manifold $\overline{C(V)}$ constructed in Section \ref{subsec: toric submanifolds}.
We first give a description of Delzant polytopes $\Delta$ in terms of \textit{manifolds with corners}.
Similar to what we have done in Section \ref{subsec: toric submanifolds}, we construct $\overline{D(V)} \subset \Delta$ from an affine subspace $V$.
In Section \ref{subsec: legendre transformation of affine linear subspaces}, we introduce the notion of a \textit{submanifold with corners}.
Similar to the case of $\overline{C(V)}$, $\overline{D(V)}$ may not be a submanifold with corners for some $V$.
We then show that $\overline{C(V)}$ is a torus-equivariantly embedded toric manifold if and only if $\overline{D(V)}$ is a submanifold with corners.

In Section \ref{subsec: torus fibration for toric manifolds}, we review how to construct a torus bundle for toric manifolds and Legendre transformation.
In Section \ref{subsec: legendre transformation of affine linear subspaces}, we construct submanifolds $\overline{D(V)}$ with corners in a Delzant polytope $\Delta$ associated to affine subspaces $V$.
In Section \ref{subsec: lyz correspondence}, we give the correspondence between $\overline{C(V)}$ and $\overline{D(V)}$.

\subsection{Torus Fibrations for Toric Manifolds}
\label{subsec: torus fibration for toric manifolds}

We review how to construct a structure of the total space of torus bundle for a toric manifold.
The base space of it should be a tropical affine manifold, which is defined as follows:
\begin{defn}
	\label{def: tropical affine manifolds}
	Let $(U; x_{1}, \ldots, x_{n}), (V;\tilde{x}_{1}, \ldots, \tilde{x}_{n})$ be local coordinates on a real manifold $B$. 
	The coordinate transformation on $U \cap V \neq \emptyset$ is a \textit{tropical affine transformation} if the coordinate transformation is given by
	\begin{equation*}
		\tilde{x} = \alpha x + \beta ~(x = {^t[x_{1} \cdots x_{n}]}, \tilde{x} = {^t[\tilde{x}_{1} \cdots \tilde{x}_{n}]})
	\end{equation*}
	with some $\alpha \in GL(n;\mathbb{Z})$ and $\beta \in \mathbb{R}^{n}$.
	A pair $(B,\{(U_{i},\varphi_{i})\}_{i \in I})$ is a tropical affine manifold if $\{(U_{i},\varphi_{i})\}_{i \in I}$ is a coordinate chart on the manifold $B$ with a tropical affine transformation.
\end{defn}

From the construction of the inhomogeneous coordinates on a toric manifold $X$, 
we can obtain a torus bundle, whose base space is a tropical affine manifold.
First, we define a tropical affine manifold $B$ from the data of a Delzant polytope $\Delta$ of $X$.

Let $B = \mathbb{R}^{n}$ and $(x_{1}, \ldots, x_{n})$ a coordinate on $B$. 
Define an open cover $\{W_{\lambda}\}_{\lambda \in \Lambda}$ of $B$ by $W_{\lambda} = \mathbb{R}^{n} = B$ for any $\lambda \in \Lambda$.
We define a homeomorphism $\psi_{\lambda}: W_{\lambda} \to \mathbb{R}^{n}$ to the image by 
\begin{equation*}
	\psi_{\lambda}(x_{1},\ldots, x_{n})
	= \left( \sum_{j=1}^{n} Q^{\lambda}_{j1}x_{j}, \ldots, \sum_{j=1}^{n} Q^{\lambda}_{jn}x_{j} \right).
\end{equation*}
Let $(x^{\lambda}_{1}, \ldots, x^{\lambda}_{n})$ be a coordinate on $\psi_{\lambda}(W_{\lambda})$ for $\lambda \in \Lambda$.
By the definition of $\psi_{\lambda}$, we have 
\begin{equation}
	(\psi_{\mu}\circ (\psi_{\lambda})^{-1})(x^{\lambda}_{1},\ldots, x^{\lambda}_{n})
	= \left(\sum_{j=1}^{n}d^{\lambda\mu}_{j1}x^{\lambda}_{j},\ldots, \sum_{j=1}^{n}d^{\lambda\mu}_{jn}x^{\lambda}_{j} \right)
\end{equation}
for all $\lambda,\mu \in \Lambda$.
Since we see in Section \ref{subsec: toric} that $D^{\lambda\mu} \in SL(n;\mathbb{Z}) \subset GL(n;\mathbb{Z})$, we see that the pair $(B,\{(W_{\lambda},\psi_{\lambda})\}_{\lambda \in \Lambda})$ is a tropical affine manifold. Hereafter, we simply write a tropical affine manifold $B$ instead of $(B,\{(W_{\lambda},\psi_{\lambda})\}_{\lambda \in \Lambda})$.

The next proposition is well-known.
\begin{prop}
	\label{prop: torus fibration from toric to tropical affine}
	The toric divisor complements $\check{M}$ of a toric manifold $X$ is a torus bundle over the tropical affine manifold $B$.
\end{prop}

The tropical affine manifold $B$ and the interior $\mathring{\Delta}$ of a Delzant polytope $\Delta$ of $X$ are diffeomorphic via the Legendre transformation.
Before giving the explanation of the Legendre transformation, we write the Guillemin potential \cite{MR1293656}.
Given a toric manifold $X$, a Delzant polytope $\Delta$ can be given by the primitive inward pointing normal vectors $u_{1}, \ldots, u_{d}$ for each facet in $\Delta$.
We write a polytope $\Delta$ explicitly as 
\begin{equation}
	\label{eq: delzant polytope}
	\Delta = \{ \xi \in (\mathfrak{t}^{n})^{\ast} \cong \mathbb{R}^{n} \mid L_{j}(\xi) \geq 0 ~\text{for}~ j = 1,\ldots, d \},
\end{equation}
where $L_{j}(\xi) = \langle \xi, u_{j} \rangle - \kappa_{j}$ for some $\kappa_{j} \in \mathbb{R}$.
We define the function $G: \mathring{\Delta} \to \mathbb{R}$ by 
\begin{equation}
	\label{eq: guillemin function}
	G(\xi) = \sum_{j=1}^{d}L_{j}(\xi)\log L_{j}(\xi),
\end{equation}
which we call the \textit{Guillemin potential} of a Delzant polytope $\Delta$.
By the theory of the Legendre transformation, we give a diffeomorphism $\Phi: \mathring{\Delta} \to B$ by 
\begin{equation*}
	\Phi (\xi_{1}, \ldots, \xi_{n})
	= \left( \frac{\partial G}{\partial \xi_{1}}, \ldots, \frac{\partial G}{\partial \xi_{n}} \right).
\end{equation*}
By \cite{MR1293656}, we obtain the potential function $F: B \to \mathbb{R}$ for the inverse of $\Phi$ by 
\begin{equation}
	\label{eq: legendre transformation}
	F(x) + G(\xi) = \sum_{i=1}^{n} x_{i}\xi_{i}.
\end{equation}
By solving Equation (\ref{eq: legendre transformation}), we can give explicitly the inverse $\Phi^{-1}:B \to \mathring{\Delta}$ by 
\begin{equation*}
	\Phi^{-1}(x_{1}, \ldots, x_{n})
	= \left( \frac{\partial F}{\partial x_{1}}, \ldots, \frac{\partial F}{\partial x_{n}} \right).
\end{equation*}
Note that this function $F$ is the potential function for the torus invariant symplectic form which we mention in Section \ref{subsec: toric}.

\subsection{Construction of Submanifolds in Delzant Polytopes}
\label{subsec: legendre transformation of affine linear subspaces}

We first consider the Legendre transformation of affine subspaces $V = \mathbb{R}p_{1} + \cdots +\mathbb{R}p_{k} + a$ in $\mathbb{R}^{n} = B$.
We write the image of $V$ under $\Phi^{-1}$ as $D(V)$. We can rewrite $D(V) = \Phi^{-1}(V) \subset \mathring{\Delta}$ as 
\begin{equation}
	\label{eq: legendre transformation of affine linear subspace}
	D(V) = \left\{ \xi \in \mathring{\Delta} \; \middle| \; \left\langle \frac{\partial G}{\partial \xi} - a, q_{l} \right\rangle = 0 ~\text{for}~ l= k+1, \ldots, n \right\}.
\end{equation}
Recall that the vectors $q_{k+1}, \ldots, q_{n} \in \mathbb{Z}^{n}$ are a primitive basis of the orthogonal subspace to $V$.

Since $V$ is a submanifold in $B$ and $\Phi^{-1}:B \to \mathring{\Delta}$ is a diffeomorphism, $D(V) = \Phi^{-1}(V)$ is a $k$-dimensional submanifold in $\mathring{\Delta}$.

In the construction of submanifolds in $\Delta$, we borrow the idea from the concepts of manifolds with corners.
We will see that a Delzant polytope $\Delta$ as a manifold with corners. 

We first review the notion of manifolds with corners in order to fix the notation.
\begin{defn}
	\label{def: smooth function boundary and corners}
	Let $U$ be an open subset in $(\mathbb{R}_{\geq 0})^{n}$. 
	A continuous function $f: U \to \mathbb{R}$ is \textit{smooth} if there exist an open subset $\tilde{U}$ in $\mathbb{R}^{n}$ and a smooth function $\tilde{f}:\tilde{U} \to \mathbb{R}$ 
	such that 
	$\tilde{U}{\mid_{(\mathbb{R}_{\geq 0})^{n}}} = U$ and $\tilde{f}{\mid_{U}} = f$.
\end{defn}

For any point $p \in (\mathbb{R}_{\geq 0})^{n} \setminus (\mathbb{R}_{>0})^{n}$, 
we define the \textit{derivative} of $f:U \to \mathbb{R}$ in Definition \ref{def: smooth function boundary and corners} at the point $p$ as 
\begin{equation*}
	\frac{\partial f}{\partial \xi_{i}}(p) = \frac{\partial \tilde{f}}{\partial \xi_{i}}(p),
\end{equation*}
where $\tilde{f}:\tilde{U} \to \mathbb{R}$ is given as in Definition \ref{def: smooth function boundary and corners} and $(\mathbb{R}_{\geq 0})^{n} \subset \mathbb{R}^{n} = \{(\xi_{1}, \ldots, \xi_{n}) \}$.

\begin{defn}
	\label{def: diffeo boundary and corners}
	Let $U$, $V$ be open subsets in $(\mathbb{R}_{\geq 0})^{n}$.
	A continuous map $f:U \to V$ is a \textit{diffeomorphism} if there exist open subsets $\tilde{U}$, $\tilde{V}$ in $\mathbb{R}^{n}$ and a diffeomorphism $\tilde{f}:\tilde{U} \to \tilde{V}$
	such that 
	$\tilde{U}{\mid_{(\mathbb{R}_{\geq 0})^{n}}} = U$, 
	$\tilde{V}{\mid_{(\mathbb{R}_{\geq 0})^{n}}} = V$, and
	$\tilde{f}{\mid_{U}} = f$.
\end{defn}
\begin{defn}
	\label{def: manifold with boundary and corners}
	A Hausdorff space $M$ is a \textit{manifold with corners} if 
	$M$ has an open covering $\{U_{\lambda} \}_{\lambda \in \Lambda}$ and a homeomorphism $\varphi_{\lambda}: U_{\lambda} \to \varphi_{\lambda}(U_{\lambda}) \subset (\mathbb{R}_{\geq 0})^{n}$ for each $\lambda \in \Lambda$ such that for $\lambda, \mu \in \Lambda$ if $U_{\lambda} \cap U_{\mu} \neq \emptyset$,
	\begin{equation*}
		\varphi_{\mu} \circ \varphi_{\lambda}^{-1}: \varphi_{\lambda} (U_{\lambda} \cap U_{\mu}) \to \varphi_{\mu} (U_{\lambda} \cap U_{\mu})
	\end{equation*}
	is a diffeomorphism in the sense of Definition \ref{def: diffeo boundary and corners}.
	We call $\{(U_{\lambda},\varphi_{\lambda})\}_{\lambda \in \Lambda}$ a system of coordinate charts on $M$.
\end{defn}

Recall that for each vertex $\lambda \in \Lambda$, we have the primitive inward pointing normal vectors $u^{\lambda}_{1}, \ldots, u^{\lambda}_{n}$ to the facets meeting the vertex $\lambda$. 
Let $\xi^{\lambda} = (\xi^{\lambda}_{1}, \ldots, \xi^{\lambda}_{n}) \in (\mathfrak{t}^{n})^{\ast} \cong \mathbb{R}^{n}$ and $L^{\lambda}_{j}(\xi^{\lambda})= \langle u^{\lambda}_{j}, \xi^{\lambda} \rangle - \kappa^{\lambda}_{j}$.
We define $\Delta_{\lambda}$ and $\mathring{\Delta}_{\lambda}$ by
\begin{align*}
	\Delta_{\lambda} 
	&= \{\xi^{\lambda} \in (\mathfrak{t}^{n})^{\ast} \mid L^{\lambda}_{j}(\xi^{\lambda}) \geq 0, \; j = 1,\ldots,n \} \cong (\mathbb{R}_{\geq 0})^{n}, \\
	\mathring{\Delta}_{\lambda} 
	&= \Delta_{\lambda} \setminus \partial \Delta_{\lambda}
	= \{ \xi^{\lambda} \in (\mathfrak{t}^{n})^{\ast} \mid L^{\lambda}_{j}(\xi^{\lambda}) > 0, \; j = 1,\ldots,n \} \cong (\mathbb{R}_{> 0})^{n}.
\end{align*}
We further define the function $G^{\lambda}:\mathring{\Delta}_{\lambda} \to \mathbb{R}$ by
\begin{equation}
	\label{eq: guillemin function local}
	G^{\lambda}(\xi^{\lambda})
	= \sum_{j=1}^{n} L^{\lambda}_{j}(\xi^{\lambda}) \log L^{\lambda}_{j}(\xi^{\lambda}).
\end{equation}
We then consider the solution $\xi^{\lambda}$ of the system of the following equations;
\begin{equation}
	\label{eq: differential equations on legendre for local}
	\frac{\partial G^{\lambda}}{\partial \xi^{\lambda}_{1}} = \frac{\partial G}{\partial \xi_{1}}, \ldots, \frac{\partial G^{\lambda}}{\partial \xi^{\lambda}_{n}} = \frac{\partial G}{\partial \xi_{n}}.
\end{equation}
We solve the system of equations (\ref{eq: differential equations on legendre for local}) in terms of $\xi$.
\begin{lem}
	\label{lemma: formula for local models}
	The solution $\xi^{\lambda}$ of the system of equations (\ref{eq: differential equations on legendre for local}) is given by
	\begin{equation}
		\label{eq: formula for local models}
		\left[
		\begin{matrix}
			\xi^{\lambda}_{1} \\
			\vdots \\
			\xi^{\lambda}_{n}
		\end{matrix}
		\right]
		=
		\left[
		\begin{matrix}
			v^{\lambda}_{1} & \cdots & v^{\lambda}_{n}
		\end{matrix}
		\right]
		\left[
		\begin{matrix} \displaystyle
			\exp \left( \left\langle v^{\lambda}_{1}, \frac{\partial G}{\partial \xi} \right\rangle -1 \right) + \kappa^{\lambda}_{1} \\
			\vdots \\
			\displaystyle
			\exp \left( \left\langle v^{\lambda}_{n}, \frac{\partial G}{\partial \xi} \right\rangle -1 \right) + \kappa^{\lambda}_{n} 
		\end{matrix}
		\right]
	\end{equation}
	for $\lambda \in \Lambda$.

	\begin{proof}
		Calculating the left-hand side of each equation in (\ref{eq: differential equations on legendre for local}),
		we obtain 
		\begin{align*}
			\left[ 
			\begin{matrix} \displaystyle
				\frac{\partial G^{\lambda}}{\partial \xi^{\lambda}_{1}} \\
				\vdots \\
				\displaystyle
				\frac{\partial G^{\lambda}}{\partial \xi^{\lambda}_{n}}
			\end{matrix}
			\right]
			&=
			\left[ 
			\begin{matrix} \displaystyle
				\sum_{j=1}^{n} \langle u^{\lambda}_{j}, e_{1} \rangle (1 + \log L^{\lambda}_{j}) \\
				\vdots \\
				\displaystyle
				\sum_{j=1}^{n} \langle u^{\lambda}_{j}, e_{n} \rangle (1 + \log L^{\lambda}_{j})
			\end{matrix} 
			\right]\\
			&= 
			\left[ 
			\begin{matrix}
				\langle u^{\lambda}_{1}, e_{1} \rangle & \cdots & \langle u^{\lambda}_{n}, e_{1} \rangle \\
				\vdots & \ddots & \vdots \\
				\langle u^{\lambda}_{1}, e_{n} \rangle & \cdots & \langle u^{\lambda}_{n}, e_{n} \rangle 
			\end{matrix} 
			\right]
			\left[ 
			\begin{matrix}
				1 + \log L^{\lambda}_{1} \\
				\vdots \\
				1 + \log L^{\lambda}_{n}
			\end{matrix}
			\right] \\
			&=
			\left[ 
			\begin{matrix}
				u^{\lambda}_{1} & \cdots & u^{\lambda}_{n}
			\end{matrix}
			\right]
			\left[
			\begin{matrix}
				1 + \log L^{\lambda}_{1} \\
				\vdots \\
				1 + \log L^{\lambda}_{n}
			\end{matrix}
			\right].
		\end{align*}
		From equations in (\ref{eq: differential equations on legendre for local}), we have 
		\begin{equation*}
			\left[ 
			\begin{matrix}
				\log L^{\lambda}_{1} \\
				\vdots \\
				\log L^{\lambda}_{n} 
			\end{matrix}
			\right]
			=
			\left[ 
			\begin{matrix}
				{^t v^{\lambda}_{1}} \\
				\vdots \\
				{^t v^{\lambda}_{n}}
			\end{matrix}
			\right]
			\left[ 
			\begin{matrix}
				\displaystyle
				\frac{\partial G}{\partial \xi_{1}} \\
				\vdots \\
				\displaystyle
				\frac{\partial G}{\partial \xi_{n}}			
			\end{matrix}
			\right]
			- 
			\left[ 
			\begin{matrix}
				1 \\
				\vdots \\
				1
			\end{matrix}
			\right]
			=
			\left[ 
			\begin{matrix}
				\displaystyle
				\left\langle v^{\lambda}_{1}, \frac{\partial G}{\partial \xi} \right\rangle - 1 \\
				\vdots \\
				\displaystyle
				\left\langle v^{\lambda}_{n}, \frac{\partial G}{\partial \xi} \right\rangle -1
			\end{matrix}
			\right].
		\end{equation*}
		Since we see
		\begin{equation*}
			\left[ 
			\begin{matrix}
				L^{\lambda}_{1} \\
				\vdots \\
				L^{\lambda}_{n}
			\end{matrix}
			\right]
			= 
			\left[
			\begin{matrix}
				{^t u^{\lambda}_{1}} \\
				\vdots \\
				{^t u^{\lambda}_{n}} 
			\end{matrix}
			\right]
			\left[
			\begin{matrix}
				\xi^{\lambda}_{1} \\
				\vdots \\
				\xi^{\lambda}_{n}
			\end{matrix}
			\right]
			-
			\left[ 
			\begin{matrix}
				\kappa^{\lambda}_{1} \\
				\vdots \\
				\kappa^{\lambda}_{n}
			\end{matrix}
			\right]
		\end{equation*}
		from the definition of $L^{\lambda}_{1}, \ldots, L^{\lambda}_{n}$, we have 
		\begin{equation*}
			\left[
				\begin{matrix}
					\xi^{\lambda}_{1} \\
					\vdots \\
					\xi^{\lambda}_{n}
				\end{matrix}
			\right]
			=
			\left[
				\begin{matrix}
					v^{\lambda}_{1} & \cdots & v^{\lambda}_{n}
				\end{matrix}
			\right]
			\left( 
				\left[
				\begin{matrix}
					L^{\lambda}_{1} \\
					\vdots \\
					L^{\lambda}_{n}
				\end{matrix}
				\right]
				+
				\left[
					\begin{matrix}
						\kappa^{\lambda}_{1} \\
					\vdots \\
					\kappa^{\lambda}_{n}
					\end{matrix}
				\right]
			\right).
		\end{equation*}
		From the equations we derived above, we obtain 
		the solution (\ref{eq: formula for local models}).
	\end{proof}
\end{lem}
In other words, we define the image of the Legendre transformation for $\mathring{\Delta}_{\lambda}$ with the potential function $G^{\lambda}$ to be the same as one of the Legendre transformation for $\mathring{\Delta}$ with the potential function $G$, which is $B$.
\begin{defn}
	\label{def: global to local transformation}
	For each $\lambda \in \Lambda$, 
	we define \textit{the transformation $\Psi_{\lambda}: \mathring{\Delta} \to \mathring{\Delta}_{\lambda}$} by Equation (\ref{eq: formula for local models}).
\end{defn}

Let $d$ be the number of facets of $\Delta$.
For each $\lambda \in \Lambda$, we define the $d-n$ facets which do not contain the vertex $\lambda$ as $\mathcal{F}^{\lambda}_{n+1}, \ldots, \mathcal{F}^{\lambda}_{d}$.
We next show that $\Psi_{\lambda}: \mathring{\Delta} \to \mathring{\Delta}_{\lambda}$ can be extended $\overline{\Psi}_{\lambda}: \Delta \setminus (\bigcup_{i=n+1}^{d} \mathcal{F}^{\lambda}_{i}) \to \Delta_{\lambda}$ such that $\overline{\Psi}_{\lambda}{\mid_{\mathring{\Delta}}} = \Psi_{\lambda}$.
\begin{lem}
	\label{lemma: boundary defined}
	From Equation (\ref{eq: formula for local models}), we obtain for $i=1,\ldots, n$
	\begin{equation}
		\label{eq: formula for local models2}
		L^{\lambda}_{i}(\xi^{\lambda})
		= \exp \left(\sum_{j=n+1}^{d}\langle v^{\lambda}_{i},u^{\lambda}_{j} \rangle \right) L^{\lambda}_{i}(\xi) \prod_{j=n+1}^{d}L^{\lambda}_{j}(\xi)^{\langle v^{\lambda}_{i}, u^{\lambda}_{j} \rangle},
	\end{equation}
	where for $j = n+1,\ldots,d$, the vector $u^{\lambda}_{j}$ is the primitive inward pointing normal vectors to $\mathcal{F}^{\lambda}_{j}$ and $L^{\lambda}_{j} (\xi) = \langle u^{\lambda}_{j}, \xi \rangle - \kappa^{\lambda}_{j}$ for some $\kappa^{\lambda}_{j} \in \mathbb{R}^{n}$.

\begin{proof}
	Since the function $G$ can be written as 
	\begin{align*}
		G(\xi) 
		&= \sum_{i=1}^{d} L_{i} (\xi) \log L_{i} (\xi) \\
		&= \sum_{i=1}^{n} L^{\lambda}_{i} (\xi) \log L^{\lambda}_{i} (\xi) + \sum_{i=n+1}^{d} L^{\lambda}_{i} (\xi) \log L^{\lambda}_{i} (\xi)
	\end{align*}
	for any $\lambda \in \Lambda$, 
	we can calculate 
	\begin{align*}
		\left[
			\begin{matrix}
				\displaystyle\frac{\partial G}{\partial \xi_{1}} \\
				\vdots \\
				\displaystyle\frac{\partial G}{\partial \xi_{n}}
			\end{matrix}
		\right]
		&= 
		\left[
			\begin{matrix}
				u^{\lambda}_{1} & \cdots & u^{\lambda}_{d}
			\end{matrix}
		\right]
		\left[
			\begin{matrix}
				\log L^{\lambda}_{1} +1\\
				\vdots \\
				\log L^{\lambda}_{d}+1
			\end{matrix}
		\right].
	\end{align*}
	Since by the definition of Delzant polytopes, 
	\begin{align*}
		\left\langle v^{\lambda}_{i},\frac{\partial G}{\partial \xi} \right\rangle  
		&=
		\left[
			\begin{matrix}
				\langle v^{\lambda}_{i},u^{\lambda}_{1} \rangle 
				& \cdots 
				& \langle v^{\lambda}_{i},u^{\lambda}_{d} \rangle
			\end{matrix}
		\right]
		\left[
			\begin{matrix}
				\log L^{\lambda}_{1} +1\\
				\vdots \\
				\log L^{\lambda}_{d}+1
			\end{matrix}
		\right] \\
		&= 
		\left[
			\begin{matrix}
				\delta_{i1}
				& \cdots 
				& \delta_{in}
				&\langle v^{\lambda}_{i},u^{\lambda}_{n+1} \rangle
				& \cdots 
				& \langle v^{\lambda}_{i},u^{\lambda}_{d} \rangle
			\end{matrix}
		\right]
		\left[
			\begin{matrix}
				\log L^{\lambda}_{1} +1\\
				\vdots \\
				\log L^{\lambda}_{d}+1
			\end{matrix}
		\right] \\
		&=
			\log L^{\lambda}_{i} +1 
			+ \sum_{j=n+1}^{d}\langle v^{\lambda}_{i},u^{\lambda}_{j} \rangle \left( \log L^{\lambda}_{j} +1\right)
	\end{align*}
	for any $i= 1,\ldots,n$, we obtain 
	\begin{equation*}
		\exp \left( \left\langle v^{\lambda}_{i},\frac{\partial G}{\partial \xi} \right\rangle -1 \right) 
		=
		\exp \left( \sum_{j=n+1}^{d}\langle v^{\lambda}_{i},u^{\lambda}_{j} \rangle \right) L^{\lambda}_{i} \prod_{j=n+1}^{d}(L^{\lambda}_{j})^{\langle v^{\lambda}_{i}, u^{\lambda}_{j} \rangle}.
	\end{equation*}
	From Equation (\ref{eq: formula for local models}), we obtain 
	\begin{align*}
		L^{\lambda}_{i}(\xi^{\lambda}) 
		&= 
		\langle u^{\lambda}_{i},\xi^{\lambda} \rangle - \kappa^{\lambda}_{i} \\
		&=
		\exp \left( \left\langle v^{\lambda}_{i},\frac{\partial G}{\partial \xi} \right\rangle -1 \right)\\
		&=
		\exp \left(\sum_{j=n+1}^{d}\langle v^{\lambda}_{i},u^{\lambda}_{j} \rangle \right) L^{\lambda}_{i} \prod_{j=n+1}^{d}(L^{\lambda}_{j})^{\langle v^{\lambda}_{i}, u^{\lambda}_{j} \rangle}
	\end{align*}
	for any $i=1,\ldots,n$.
	Therefore, we obtain Equation (\ref{eq: formula for local models2}).
\end{proof}
\end{lem}

From Equation (\ref{eq: formula for local models2}),
we rewrite the expression of Equation 
(\ref{eq: formula for local models}).

\begin{lem}
	\label{lemma: solution for global to local}
	From Equation (\ref{eq: formula for local models2}), we have 
	\begin{equation}
		\label{eq: formula for local models3}
		\left[
			\begin{matrix}
				\xi^{\lambda}_{1} \\
				\vdots \\
				\xi^{\lambda}_{n}
			\end{matrix}
		\right]
		=
		Q^{\lambda}
		\left[
			\begin{matrix}
				\displaystyle
				\exp \left(\sum_{j=n+1}^{d}\langle v^{\lambda}_{1},u^{\lambda}_{j} \rangle \right) L^{\lambda}_{1}(\xi) \prod_{j=n+1}^{d}L^{\lambda}_{j}(\xi)^{\langle v^{\lambda}_{1}, u^{\lambda}_{j} \rangle} + \kappa^{\lambda}_{1} \\
				\vdots \\
				\displaystyle 
				\exp \left(\sum_{j=n+1}^{d}\langle v^{\lambda}_{n},u^{\lambda}_{j} \rangle \right) L^{\lambda}_{n}(\xi) \prod_{j=n+1}^{d}L^{\lambda}_{j}(\xi)^{\langle v^{\lambda}_{n}, u^{\lambda}_{j} \rangle} + \kappa^{\lambda}_{n}
			\end{matrix}
		\right],
	\end{equation}
	where for $j = n+1,\ldots,d$ the vector $u^{\lambda}_{j}$ is the inward pointing normal vectors to $\mathcal{F}^{\lambda}_{j}$ and $L^{\lambda}_{j} (\xi) = \langle u^{\lambda}_{j}, \xi \rangle - \kappa^{\lambda}_{j}$ for some $\kappa^{\lambda}_{j} \in \mathbb{R}^{n}$.

\begin{proof}
	Recall that $Q^{\lambda} = \left[
		\begin{matrix}
			v^{\lambda}_{1} & \cdots &v^{\lambda}_{n}
		\end{matrix}
	\right]$, where $v^{\lambda}_{1}, \ldots, v^{\lambda}_{n}$ are the direction vectors of the edges from the vertex $\lambda$.

	Since $L^{\lambda}_{i}(\xi^{\lambda}) = \langle u^{\lambda}_{i},\xi^{\lambda} \rangle - \kappa^{\lambda}_{i}$ for each $i = 1,\ldots, n$, we have 
	\begin{align*}
		\left[
			\begin{matrix}
				L^{\lambda}_{1}(\xi^{\lambda}) \\
				\vdots \\
				L^{\lambda}_{n}(\xi^{\lambda})
			\end{matrix}
		\right]
		&=
		\left[
			\begin{matrix}
				\langle u^{\lambda}_{1},\xi^{\lambda} \rangle - \kappa^{\lambda}_{1} \\
				\vdots \\
				\langle u^{\lambda}_{n},\xi^{\lambda} \rangle - \kappa^{\lambda}_{n}
			\end{matrix}
		\right] \\
		&=
		\left[
			\begin{matrix}
				{^t u^{\lambda}_{1}} \\
				\vdots \\
				{^t u^{\lambda}_{n}}
			\end{matrix}
		\right]
		\left[
			\begin{matrix}
				\xi^{\lambda}_{1} \\
				\vdots \\
				\xi^{\lambda}_{n}
			\end{matrix}
		\right]
		- 
		\left[
			\begin{matrix}
				\kappa^{\lambda}_{1} \\
				\vdots \\
				\kappa^{\lambda}_{n}
			\end{matrix}
		\right].
	\end{align*}
	Since the inverse matrix of the square matrix ${^t[ \begin{matrix} u^{\lambda}_{1} &\cdots& u^{\lambda}_{n} \end{matrix}]}$ is $Q^{\lambda}$ (see Lemma \ref{lemma: direction vectors vs normal vectors}), we obtain 
	\begin{equation*}
		\left[
			\begin{matrix}
				\xi^{\lambda}_{1} \\
				\vdots \\
				\xi^{\lambda}_{n}
			\end{matrix}
		\right]
		=
		Q^{\lambda}
		\left(\left[
			\begin{matrix}
				L^{\lambda}_{1}(\xi^{\lambda}) \\
				\vdots \\
				L^{\lambda}_{n}(\xi^{\lambda})
			\end{matrix}
		\right]
		+
		\left[
			\begin{matrix}
				\kappa^{\lambda}_{1} \\
				\vdots \\
				\kappa^{\lambda}_{n}
			\end{matrix}
		\right]\right).
	\end{equation*}
	From Equation (\ref{eq: formula for local models2}), we see Equation (\ref{eq: formula for local models3}) holds.
\end{proof}
\end{lem}

We define the transformation $\Psi_{\lambda}: \mathring{\Delta} \to \mathring{\Delta}_{\lambda}$ by Equation (\ref{eq: formula for local models}) for each $\lambda \in \Lambda$.
Let $\mathcal{F}^{\lambda}_{1},\ldots, \mathcal{F}^{\lambda}_{n}$ be facets meeting the vertex $\lambda$.
Then, we see that $\xi \in \mathcal{F}^{\lambda}_{i}$ is equivalent to the condition that $L^{\lambda}_{i}(\xi) = 0$.
From Equation (\ref{eq: formula for local models2}), we see that $L^{\lambda}_{i}(\xi^{\lambda}) = 0$ is equivalent to $L^{\lambda}_{i}(\xi) = 0$ for any point $\xi$ in the boundary $\partial (\Delta \setminus (\bigcup_{j=n+1}^{d}\mathcal{F}^{\lambda}_{j}))$ of $\Delta \setminus (\bigcup_{j=n+1}^{d}\mathcal{F}^{\lambda}_{j})$ and each $i=1,\ldots, n$.
Thus for any $\xi \in \partial (\Delta \setminus (\bigcup_{j=n+1}^{d}\mathcal{F}^{\lambda}_{j}))$, $\xi^{\lambda} \in \partial \Delta_{\lambda}$ is defined via Equation (\ref{eq: formula for local models3}).
As a conclusion, 
using Lemma \ref{lemma: boundary defined} and Lemma \ref{lemma: solution for global to local}, we can define the map $\overline{\Psi}_{\lambda}: \Delta \setminus (\bigcup_{i=n+1}^{d} \mathcal{F}^{\lambda}_{i}) \to \Delta_{\lambda}$ by the extension of $\Psi_{\lambda}:\mathring{\Delta} \to \mathring{\Delta}_{\lambda}$.
\begin{defn}
	\label{def: global to local with boundary}
	For each $\lambda \in \Lambda$, 
	we define \textit{the transformation $\overline{\Psi}_{\lambda}:\Delta \setminus (\bigcup_{i=n+1}^{d} \mathcal{F}^{\lambda}_{i}) \to \Delta_{\lambda}$} by Equation (\ref{eq: formula for local models3}).
\end{defn}

Since $\Delta_{\lambda}$ is diffeomorphic to $(\mathbb{R}_{\geq 0})^{n}$ for each $\lambda \in \Lambda$,
we regard $\Delta$ as a manifold with corners, by using $\{(\overline{\Psi}_{\lambda}^{-1}(\Delta_{\lambda}), \overline{\Psi}_{\lambda}) \}_{\lambda \in \Lambda}$.

As a preparation for our main result,
we introduce the notion of \textit{submanifolds with corners} in manifolds with corners.
\begin{defn}
	\label{def: submanifolds with boundary and corners}
	Let $M$ be an $n$-dimensional manifold with corners, $\{(U_{\lambda},\varphi_{\lambda}) \}_{\lambda \in \Lambda}$ a system of coordinate charts on $M$.  
	A subset $N$ in $M$ is a \textit{submanifold with corners} if 
	for each $\lambda \in \Lambda$, there exist functions $g^{\lambda}_{1}, \ldots, g^{\lambda}_{n-k}$ which are all smooth in the sense of Definition \ref{def: smooth function boundary and corners} on $U_{\lambda}$
	such that 
	\begin{equation*}
		N \cap U_{\lambda} = \{p \in U_{\lambda} \mid g^{\lambda}_{1}(p) = \cdots = g^{\lambda}_{n-k}(p) = 0 \}
	\end{equation*}
	and the rank of the Jacobian matrix for $[g^{\lambda}_{1}\; g^{\lambda}_{2}\; \cdots \; g^{\lambda}_{n-k}]$ is equal to $n-k$ at any point in $N\cap U_{\lambda}$.
\end{defn}

We will observe the closure $\overline{D(V)}$ of $D(V)$ to $\Delta$ as a submanifold with corners. 
We define $D_{\lambda}(V)$ as 
\begin{equation}
	D_{\lambda}(V)
	= \left\{ \xi^{\lambda} \in \mathring{\Delta}_{\lambda} \; \middle| \; \left\langle \frac{\partial G^{\lambda}}{\partial \xi^{\lambda}} - a, q_{l} \right\rangle = 0 ~\text{for}~ l= k+1, \ldots, n \right\}.
\end{equation} 
From Equation (\ref{eq: guillemin function local}), we have 
\begin{equation*}
	\frac{\partial G^{\lambda}}{\partial \xi^{\lambda}_{i}}
	= 
	\sum_{j=1}^{n} \left(
		\langle u^{\lambda}_{j},e_{i} \rangle \log L^{\lambda}_{j} + \langle u^{\lambda}_{j},e_{i} \rangle
	\right)
\end{equation*}
for $i=1,\ldots,n$. 
Direct calculation gives us 
\begin{align*}
	\left\langle \frac{\partial G^{\lambda}}{\partial \xi^{\lambda}}, q_{l} \right\rangle 
	=& \sum_{i=1}^{n} \frac{\partial G^{\lambda}}{\partial \xi^{\lambda}_{i}} \langle e_{i},q_{l} \rangle \\
	=& \sum_{i=1}^{n}\sum_{j=1}^{n} \left(\log L^{\lambda}_{j} + 1 \right) \langle u^{\lambda}_{j},e_{i} \rangle \langle e_{i},q_{l} \rangle\\
	=& \sum_{j=1}^{n} \left(\log L^{\lambda}_{j} + 1 \right) \langle u^{\lambda}_{j},q_{l} \rangle\\
	=& \log \left( \prod_{j=1}^{n} (eL^{\lambda}_{j})^{\langle u^{\lambda}_{j}, q_{l} \rangle} \right).
\end{align*}
From this calculation, the defining equations for $D_{\lambda}(V)$ can be written as 
\begin{equation*}
	\log \left( \prod_{j=1}^{n} (eL^{\lambda}_{j})^{\langle u^{\lambda}_{j}, q_{l} \rangle} \right) 
	-
	\langle a, q_{l} \rangle
	= 0
\end{equation*}
for $l = k+1,\ldots, n$.
Hence, $D_{\lambda}(V)$ is a zero locus of 
\begin{equation*}
	\prod_{j=1}^{n} (eL^{\lambda}_{j})^{\langle u^{\lambda}_{j}, q_{l} \rangle} 
	-
	e^{\langle a, q_{l} \rangle}
	=
	0
\end{equation*}
for $l = k+1,\ldots, n$.
Recall the definitions of $\mathcal{I}^{+}_{\lambda,l}$ and $\mathcal{I}^{-}_{\lambda,l}$ given in (\ref{eq: plus}) and (\ref{eq: minus}) respectively. Then we can define the closure $\overline{D_{\lambda}(V)}$ by
\begin{equation}
	\overline{D_{\lambda}(V)}
	= \{\xi^{\lambda} \in \Delta_{\lambda} \mid g^{\lambda}_{l} (\xi^{\lambda}) = 0,\; l = k+1,\ldots, n \}
\end{equation}
where $g^{\lambda}_{l}$ is given by
\begin{equation}
	\label{eq: defining eqautions for symplectic side}
	g^{\lambda}_{l}
	= \prod_{i \in \mathcal{I}^{+}_{\lambda,l}} (eL^{\lambda}_{i})^{\langle u^{\lambda}_{i}, q_{l} \rangle}
	- e^{\langle a, q_{l} \rangle} \prod_{i \in \mathcal{I}^{-}_{\lambda,l}} (eL^{\lambda}_{i})^{-\langle u^{\lambda}_{i},q_{l} \rangle}
\end{equation}
for each $l = k+1, \ldots, n$.
Thus we obtain $\overline{D(V)} = \bigcup_{\lambda \in \Lambda} \overline{\Psi}^{-1}_{\lambda}(\overline{D_{\lambda}(V)}) \subset \Delta$.
On $\mathring{\Delta}_{\lambda}$, the defining equations of $D_{\lambda}(V)$ and $\overline{D_{\lambda}(V)}$ are equivalent, but since $L^{\lambda}_{i} = 0$ on $\mathcal{F}^{\lambda}_{i}$, we use the equations $g^{\lambda}_{l} =0$ instead of the equations $\prod_{j=1}^{n}(eL^{\lambda}_{j})^{\langle u^{\lambda}_{j},q_{l} \rangle} - e^{\langle a,q_{l} \rangle} = 0$ to define $\overline{D_{\lambda}(V)}$.

Note that by the definition of submanifolds with corners, $\overline{D(V)}$ is a $k$-dimensional submanifold with corners in $\Delta$ if the rank of the Jacobian matrix of $g^{\lambda}:= [g^{\lambda}_{k+1}\; \cdots\; g^{\lambda}_{n}]$ is equal to $n-k$ for any points on $\overline{D_{\lambda}(V)}$ and any $\lambda \in \Lambda$.

\subsection{Correspondence between Torus-equivariantly Embedded Toric Manifolds and Submanifolds in Delzant Polytopes}
\label{subsec: lyz correspondence}

We show the main result in this paper:
\begin{thm}
\label{thm: lyz correspondence}
	The rank of the Jacobian matrix $Df^{\lambda}$ is equal to the rank of the Jacobian matrix $Dg^{\lambda}$ for any $\lambda \in \Lambda$.
	In particular,
	$\overline{C(V)}$ is a complex $k$-dimensional submanifold in $X$ if and only if $\overline{D(V)}$ is a $k$-dimensional submanifold with corners in $\Delta$.
\end{thm}
$f^{\lambda} =[f^{\lambda}_{k+1} \; f^{\lambda}_{k+2}\; \cdots \; f^{\lambda}_{n}]$ and $g^{\lambda} = [g^{\lambda}_{k+1}\; g^{\lambda}_{k+2}\; \cdots \; g^{\lambda}_{n}]$ are given in Definition \ref{def: compactification} and Equation (\ref{eq: defining eqautions for symplectic side}), respectively.

We give a lemma before giving the proof of this theorem.
\begin{lem}
	\label{lemma: correspondence points boundary and zeros}
	In our case, $\vert z^{\lambda}_{i} \vert = L^{\lambda}_{i}(\xi^{\lambda})$ for any $i= 1,\ldots,n$.
	In particular, $z^{\lambda}_{i} = 0$ is equivalent to $L^{\lambda}_{i}(\xi^{\lambda}) = 0$ for any $i = 1,\ldots, n$ and $\lambda \in \Lambda$.

\begin{proof}
	In Section \ref{subsec: torus fibration for toric manifolds}, 
	we give the Legendre transformation by 
	\begin{equation*}
		x_{1} = \frac{\partial G}{\partial \xi_{1}}, 
		\ldots, 
		x_{n} = \frac{\partial G}{\partial \xi_{n}},
	\end{equation*}
	where the function $G$ is the Guillemin potential determined by Equation (\ref{eq: guillemin function}).
	From Equation (\ref{eq: differential equations on legendre for local}),
	we further calculate 
	\begin{align*}
		x_{j} 
		&=
		\frac{\partial G}{\partial \xi_{j}} \\
		&=
		\frac{\partial G^{\lambda}}{\partial \xi^{\lambda}_{j}} \\
		&=
		\sum_{l = 1}^{n}\langle u^{\lambda}_{l},e_{j} \rangle (1 + \log L^{\lambda}_{l}) \\
		&=
		\log \left(\prod_{l=1}^{n}(e L^{\lambda}_{l})^{\langle u^{\lambda}_{l},e_{j} \rangle} \right),
	\end{align*}
	for each $j = 1,\ldots, n$.
	Note that $L^{\lambda}_{1}(\xi^{\lambda}) \geq 0, \ldots, L^{\lambda}_{n}(\xi^{\lambda}) \geq 0$ for any $\xi^{\lambda} \in \Delta_{\lambda}$.
	Since $z_{j} = \exp(x_{j} + \sqrt{-1}y_{j})$ for $j = 1,\ldots, n$, we see from \cite[Definition 3.16]{yamaguchi2023toric} that 
	\begin{align*}
		z^{\lambda}_{i}
		&= 
		\prod_{j=1}^{n}(z_{j})^{Q^{\lambda}_{ji}} \\
		&=
		\exp \left(\sum_{j=1}^{n} \left(Q^{\lambda}_{ji} x_{j} + \sqrt{-1}Q^{\lambda}_{ji} y_{j} \right) \right).
	\end{align*}
	Thus we can calculate $\vert z^{\lambda}_{i} \vert$ as
	\begin{align*}
		\vert z^{\lambda}_{i} \vert
		&=
		\exp \left(\sum_{j=1}^{n}Q^{\lambda}_{ji}x_{j} \right) \\
		&=
		\exp \left(\sum_{j=1}^{n} Q^{\lambda}_{ji} \log \left( \prod_{l=1}^{n}(eL^{\lambda}_{l})^{\langle u^{\lambda}_{l},e_{j} \rangle}\right) \right) \\
		&=
		\exp \left(\log \prod_{j=1}^{n}\prod_{l=1}^{n}(eL^{\lambda}_{l})^{\langle u^{\lambda}_{l},e_{j} \rangle Q^{\lambda}_{ji}} \right) \\
		&=
		\prod_{j=1}^{n}\prod_{l=1}^{n}(eL^{\lambda}_{l})^{\langle u^{\lambda}_{l},e_{j} \rangle Q^{\lambda}_{ji}} \\
		&=
		\prod_{l=1}^{d}(eL_{l})^{\langle u^{\lambda}_{l}, v^{\lambda}_{i} \rangle} \\
		&=
		eL^{\lambda}_{i}.
	\end{align*}
	Note that since we define the matrix $Q^{\lambda} = [\begin{matrix} v^{\lambda}_{1} & \cdots & v^{\lambda}_{n} \end{matrix}]$, we see that $Q^{\lambda}_{ji} = \langle e_{j},v^{\lambda}_{i} \rangle$ for any $i,j = 1,\ldots, n$, and that from Lemma \ref{lemma: direction vectors vs normal vectors}, we have $\langle u^{\lambda}_{i}, v^{\lambda}_{j} \rangle = \delta_{ij}$ for any $i,j = 1,\ldots, n$.

	Therefore, we finally obtain
	\begin{equation*}
		\vert z^{\lambda}_{1} \vert = eL^{\lambda}_{1}(\xi^{\lambda}), 
		\ldots, 
		\vert z^{\lambda}_{n} \vert = eL^{\lambda}_{n}(\xi^{\lambda}),
	\end{equation*}
	for any $\lambda \in \Lambda$.
	In particular, we see that $z^{\lambda}_{i} = 0$ is equivalent to $L^{\lambda}_{i}(\xi^{\lambda}) = 0$ for any $i= 1,\ldots, n$ and $\lambda \in \Lambda$. 
\end{proof}
\end{lem}

Now we prove Theorem \ref{thm: lyz correspondence}.

\begin{proof}(of Theorem \ref{thm: lyz correspondence})

	By Remark \ref{remark: toric submanifold} and the definition of submanifolds with corners, we show that the two matrices $Df^{\lambda}$ and $Dg^{\lambda}$ have the same rank for any $\lambda \in \Lambda$.
	
	We define $C^{\lambda}_{ij}$ by 
	\begin{equation*}
		C^{\lambda}_{ij} =
		\begin{cases*}
			0 & if $j \in \mathcal{I}^{+}_{\lambda,i}$, \\
			1 & if $j \in \mathcal{I}^{-}_{\lambda,i} \setminus \mathcal{I}^{0}_{\lambda,i}$,
		\end{cases*}
	\end{equation*}
	for any $i = k+1,\ldots,n $, $j = 1,\ldots, n$, and $\lambda \in \Lambda$.

	Since $f^{\lambda}_{k+1}, \ldots, f^{\lambda}_{n}$ are defined in Equation (\ref{eq: defining eqautions for toric submanifolds}), we calculate 
	\begin{equation*}
		\frac{\partial f^{\lambda}_{i}}{\partial z^{\lambda}_{j}}
		= 
		\langle u^{\lambda}_{j}, q_{i} \rangle 
		e^{\langle a,q_{i} \rangle C^{\lambda}_{ij}}
		\prod_{l \in \mathcal{I}^{\pm}_{\lambda,i}} (z^{\lambda}_{l})^{\pm \langle u^{\lambda}_{l},q_{i} \rangle - \delta_{lj}}
	\end{equation*}
	for $i = k+1, \ldots, n$ and $j = 1, \ldots,n$. 
	We obtain the Jacobian matrix $Df^{\lambda}$ of $f^{\lambda}= [f^{\lambda}_{k+1}\; \cdots\; f^{\lambda}_{n}]$ is as follows:
	\begin{equation}
		\label{eq: Jacobian matrix of complex side}
		Df^{\lambda}(z^{\lambda})
		= \left[
			\langle u^{\lambda}_{j}, q_{i} \rangle 
			e^{\langle a,q_{i} \rangle C^{\lambda}_{ij}}
			\prod_{l \in \mathcal{I}^{\pm}_{\lambda,i}} (z^{\lambda}_{l})^{\pm \langle u^{\lambda}_{l},q_{i} \rangle - \delta_{lj}}
		\right].
	\end{equation}

	For the Jacobian matrix $Dg^{\lambda}$, since $g^{\lambda}_{k+1}, \ldots, g^{\lambda}_{n}$ are defined in (\ref{eq: defining eqautions for symplectic side}), 
	we calculate
	\begin{align*}
		\frac{\partial g^{\lambda}_{i}}{\partial \xi^{\lambda}_{j}}
		=& 
		\sum_{m=1}^{n}\prod_{l \in \mathcal{I}^{+}_{\lambda,i}} \langle u^{\lambda}_{m}, q_{i} \rangle \langle u^{\lambda}_{m}, e_{j} \rangle e^{\langle u^{\lambda}_{l},q_{i} \rangle} (L^{\lambda}_{l})^{\langle u^{\lambda}_{l}, q_{i} \rangle - \delta_{lm}} \\
		&+ e^{\langle a,q_{i} \rangle C^{\lambda}_{ij}} \sum_{m=1}^{n}\prod_{l \in \mathcal{I}^{-}_{\lambda,i}} \langle u^{\lambda}_{m}, q_{i} \rangle \langle u^{\lambda}_{m}, e_{j} \rangle e^{- \langle u^{\lambda}_{l},q_{i} \rangle} (L^{\lambda}_{l})^{- \langle u^{\lambda}_{l}, q_{i} \rangle - \delta_{lm}}
	\end{align*}
	for $l = k+1, \ldots, n$ and $j = 1, \ldots, n$.
	Since the Jacobian matrix $Dg^{\lambda}$ of $g^{\lambda}=[g^{\lambda}_{j+1}\; \cdots \; g^{\lambda}_{n}]$ is expressed as 
	\begin{align*}
		Dg^{\lambda}
		&=
		\left[
			\begin{matrix}
				\displaystyle
				\frac{\partial g^{\lambda}_{k+1}}{\partial \xi^{\lambda}_{1}} 
				& \cdots & 
				\displaystyle \frac{\partial g^{\lambda}_{k+1}}{\partial \xi^{\lambda}_{n}} \\
				\vdots & \ddots & \vdots \\
				\displaystyle
				\frac{\partial g^{\lambda}_{n}}{\partial \xi^{\lambda}_{1}} 
				& \cdots & 
				\displaystyle
				\frac{\partial g^{\lambda}_{n}}{\partial \xi^{\lambda}_{n}}
			\end{matrix}
		\right] \\
		&=
		\left[
			\langle u^{\lambda}_{j}, q_{i} \rangle
			e^{\langle a,q_{i} \rangle C^{\lambda}_{ij}}
			\prod_{l \in \mathcal{I}^{\pm}_{\lambda,i}} e^{\pm \langle u^{\lambda}_{l}, q_{i} \rangle}(L^{\lambda}_{l})^{\pm \langle u^{\lambda}_{l}, q_{i} \rangle - \delta_{lj}} 
		\right]
		\left[
			\begin{array}{c}
				^t u^{\lambda}_{1} \\
				^t u^{\lambda}_{2} \\
				\vdots \\
				^t u^{\lambda}_{n}
			\end{array}
		\right],
	\end{align*}
	we obtain the Jacobian matrix $Dg^{\lambda}$ of $g^{\lambda}= [g^{\lambda}_{k+1}\; \cdots\; g^{\lambda}_{n}]$ as follows;
	\begin{equation}
	\label{eq: Jacobian matrix of symplectic side}
		Dg^{\lambda}
		= e \left[
			\langle u^{\lambda}_{j}, q_{i} \rangle
			e^{\langle a,q_{i} \rangle C^{\lambda}_{ij}}
			\prod_{l \in \mathcal{I}^{\pm}_{\lambda,i}} (eL^{\lambda}_{l})^{\pm \langle u^{\lambda}_{l}, q_{i} \rangle - \delta_{lj}} 
		\right]
		\left[
			\begin{array}{c}
				^t u^{\lambda}_{1} \\
				^t u^{\lambda}_{2} \\
				\vdots \\
				^t u^{\lambda}_{n}
			\end{array}
		\right].
	\end{equation}
	
	By the definition of Delzant polytopes, the $n \times n$ matrix $^t[u^{\lambda}_{1} \cdots u^{\lambda}_{n}]$ does not influence the rank of the Jacobian matrix $Dg^{\lambda}$. 
	Since Lemma \ref{lemma: correspondence points boundary and zeros} shows that $\vert z^{\lambda}_{i} \vert  = eL^{\lambda}_{i} (\xi^{\lambda})$ for $i = 1,\ldots, n$,
	the rank of the matrix $Df^{\lambda}(\vert z^{\lambda} \vert)$ is equal to the rank of the matrix $Dg^{\lambda}$.
	Corollary \ref{cor: subtorus action on Jacobian matrix2} shows that the rank of the matrix $Df^{\lambda}(\vert z^{\lambda} \vert)$ is equal to the rank of the matrix $Df^{\lambda}(z^{\lambda})$ for any point $z^{\lambda}$.
	Therefore, comparing Equation (\ref{eq: Jacobian matrix of complex side}) with Equation (\ref{eq: Jacobian matrix of symplectic side}), 
	we prove that the rank of $Df^{\lambda}$ is equal to that of $Dg^{\lambda}$ for any $\lambda \in \Lambda$.
\end{proof}

We summarize what we have done.
In Section \ref{subsec: torus fibration for toric manifolds}, we review the construction of a torus bundle for the toric divisor complements $\check{M}$ of a toric manifold $X$ over a tropical affine manfiold $B$ and the Legendre transformation with respect to the Guillemin potential on the interior $\mathring{\Delta}$ of a Delzant polytope of $X$.
Actually, the torus fibrations $\check{\pi}: \check{M} \to B$, the diffeomorphisms $\Phi$ and the moment map $\mu$ for the $T^{n}$-action on $\check{M} \subset X$ are commutative:
\begin{equation*}
	\begin{tikzcd}
		\check{M} \arrow[rd, "\check{\pi}"'] \arrow[rr, "\mu|_{\check{M}}"] & & \mathring{\Delta} \arrow[ld, "\Phi"] \\
		& B &
	\end{tikzcd}
\end{equation*}
from the theory of toric manifolds.
In Section \ref{subsec: legendre transformation of affine linear subspaces}, we construct a submanifold $D(V) = \Phi^{-1}(V)$ in $\mathring{\Delta}$ by using the diffeomorphism $\Phi$ for an affine subspace $V$ in $\mathbb{R}^{n} \cong B$. 
Comparing with the construction of complex submanifolds $C(V)$ given in Section \ref{subsec: toric submanifolds},
we have 
\begin{equation*}
	\begin{tikzcd}
		C(V) \arrow[rd, "\check{\pi}|_{C(V)}"'] \arrow[rr, "\mu|_{C(V)}"] & & D(V) \arrow[ld, "\Phi|_{D(V)}"] \\
		& V &
	\end{tikzcd}
\end{equation*}
from the theory of toric manifolds.
Theorem \ref{thm: lyz correspondence} shows that under the following description:
\begin{equation*}
	\begin{tikzcd}
		\overline{C(V)} \arrow[r,hook] & X \arrow[rr,"\mu"] & & \Delta & \overline{D(V)} \arrow[l, hook] \\
		C(V) \arrow[r, hook] \arrow[rrdd, "\check{\pi}|_{C(V)}"'] \arrow[u, hook] & \check{M} \arrow[rd, "\check{\pi}"'] \arrow[rr, "\mu|_{\check{M}}"] \arrow[u, hook] & &\mathring{\Delta} \arrow[ld, "\Phi"] \arrow[u, hook] & D(V) \arrow[lldd, "\Phi|_{D(V)}"] \arrow[l, hook] \arrow[u, hook] \\
		& & B & & \\
		& & V \arrow[u, hook] & &
	\end{tikzcd}	
\end{equation*}
the conditions that $\overline{C(V)}$ is a complex submanifold in $X$ is equivalent to the conditions that $\overline{D(V)}$ is a submanifold with corners.

\section{Standard Local Models to Classify Submanifolds in Delzant Polytopes}
\label{subsec: classification of submanifolds in delzant}

We give a \textit{standard local model} $\mathcal{V}_{O}$ for our construction of $\overline{D(V)}$ 
for us to classify more easily the conditions that $\overline{D(V)}$ is a submanifold with corners.

Define $\Delta_{O}$ by
\begin{equation*}
	\Delta_{O}
	= \{ \xi^{O} \in (\mathfrak{t}^{n})^{\ast} \mid \xi^{O}_{1} \geq 0, \ldots, \xi^{O}_{n} \geq 0 \} \cong (\mathbb{R}_{\geq 0})^{n}.
\end{equation*}
For a $k$-dimensional affine subspace $V = \mathbb{R}p_{1} + \cdots + \mathbb{R}p_{k} + a$ in $\mathbb{R}^{n}$ with a primitive basis $p_{1}, \ldots, p_{k} \in \mathbb{Z}^{n}$, 
we choose a primitive basis $q_{k+1}, \ldots, q_{n} \in \mathbb{Z}^{n}$ of the orthogonal subspace to $V$.
We can write $\overline{D_{O}(V)}$ as 
\begin{equation*}
	\overline{D_{O}(V)}
	= \{\xi^{O} \in \Delta_{O} \mid g^{O}_{j} (\xi^{O}) = 0,\; j=k+1,\ldots,n \},
\end{equation*}
where the functions $g^{O}_{k+1},\ldots, g^{O}_{n}$ are defined by 
\begin{equation*}
	g^{O}_{j}(\xi^{O})
	= 
	\prod_{i \in \mathcal{I}^{+}_{O,j}} (e \xi^{O}_{i})^{\langle e_{i},q_{j} \rangle} 
	- e^{\langle a,q_{j}\rangle} \prod_{i \in \mathcal{I}^{-}_{O,j}} (e \xi^{O}_{i})^{-\langle e_{i},q_{j} \rangle}
\end{equation*}
for $j=k+1,\ldots,n$. Here, the sets $\mathcal{I}^{+}_{O,j}$, $\mathcal{I}^{-}_{O,j}$ are defined by 
\begin{align}
	\label{eq: plus local}
	&\mathcal{I}^{+}_{O,j} = 
	\{i \in \{1,\ldots, n \} \mid \langle e_{i},q_{j} \rangle \geq 0 \}, \\
	\label{eq: minus local}
	&\mathcal{I}^{-}_{O,j} = 
	\{i \in \{1,\ldots, n \} \mid \langle e_{i},q_{j} \rangle \leq 0 \}. 
\end{align}

\begin{defn}
	\label{def: trivial local model}
	Fix $n \geq 2$ and $k = 1,\ldots, n-1$.
	\textit{The standard local model $\mathcal{V}_{O}$ for $n$ and $k$} is defined to be the set of the $k$-dimensional affine subspaces $V$ such that the rank of the Jacobian matrix $Dg^{O}$ of $g^{O}:= [g^{O}_{k+1}\; \cdots \; g^{O}_{n}]$ is equal to $n-k$.
\end{defn}

By using the standard local models, one can check whether $\overline{D(V)}$ is a submanifold with corners for a given affine subspace $V$ in the following manner.

Since each vertex $\lambda \in \Lambda$ of a Delzant polytope of a compact toric manifold has direction vectors which form a basis of $\mathbb{Z}^{n}$,
we can translate the direction vectors $v^{\lambda}_{1},\ldots,v^{\lambda}_{n}$ to the standard basis $e_{1},\ldots, e_{n}$.
Thus we define the following transformation.
\begin{defn}
	\label{def: standard local to local}
	We define \textit{the transformation $\Delta_{O} \to \Delta_{\lambda}$} by 
	\begin{equation}
		\label{eq: transformation from standard local to local}
		(\xi^{O}_{1}, \ldots, \xi^{O}_{n}) 
		\mapsto 
		(Q^{\lambda}(\xi^{O}_{1} + \kappa^{\lambda}_{1}), \ldots, Q^{\lambda}(\xi^{O}_{n} + \kappa^{\lambda}_{n}))
	\end{equation}
	for each $\lambda \in \Lambda$.
\end{defn}
\begin{lem}
	\label{lemma: boundary of standard local to boundary of local}
	From the transformation given in Equation (\ref{eq: transformation from standard local to local}), we obtain 
	\begin{equation*}
		L^{\lambda}_{1}(\xi^{\lambda})
		= \xi^{O}_{1},\ldots, 
		L^{\lambda}_{n}(\xi^{\lambda})
		= \xi^{O}_{n}
	\end{equation*}
	for any $\lambda \in \Lambda$.
\begin{proof}
	From the transformation given in Equation (\ref{eq: transformation from standard local to local}) and the definition of $L^{\lambda}(\xi^{\lambda})$, we calculate
	\begin{align*}
		L^{\lambda}_{i}(Q^{\lambda}(\xi^{O}_{1} + \kappa^{\lambda}_{1}), \ldots, Q^{\lambda}(\xi^{O}_{n} + \kappa^{\lambda}_{n}))
		&=
		\left\langle u^{\lambda}_{i}, \left[ \begin{matrix}
			Q^{\lambda}(\xi^{O}_{1} + \kappa^{\lambda}_{1}) \\
			\vdots \\ 
			Q^{\lambda}(\xi^{O}_{n} + \kappa^{\lambda}_{n})
		\end{matrix}\right] \right\rangle - \kappa^{\lambda}_{i} \\
		&=
		\langle u^{\lambda}_{i}, Q^{\lambda}(\xi^{O}+\kappa^{\lambda}) \rangle - \kappa^{\lambda}_{i} \\
		&=
		\langle {^t Q^{\lambda}} u^{\lambda}_{i}, \xi^{O}+\kappa^{\lambda} \rangle - \kappa^{\lambda}_{i} \\
		&=
		\left\langle \left[
			\begin{matrix}
				{^t v^{\lambda}_{1}} \\
				\vdots \\
				{^t v^{\lambda}_{n}}
			\end{matrix}
		\right] u^{\lambda}_{i}, \xi^{O}+\kappa^{\lambda} \right\rangle - \kappa^{\lambda}_{i} \\
		&=
		\langle e_{i}, \xi^{O}+\kappa^{\lambda} \rangle - \kappa^{\lambda}_{i} \\
		&= \xi^{O}_{i}.
	\end{align*}
	Theorefore, we see that $L^{\lambda}_{1} = \xi^{O}_{1},\ldots, L^{\lambda}_{n} = \xi^{O}_{n}$.
\end{proof}
\end{lem}

Since we constructed $\overline{D(V)}$ locally around each vertex, the following proposition can be verified.
\begin{prop}
	\label{prop: local models to global}
	Assume that for an affine subspace $V = \mathbb{R}p_{1}+ \cdots + \mathbb{R}p_{k} + a$ in $\mathbb{R}^{n}$, the rank of the Jacobian matrix $Dg^{O}$ is equal to $n-k$.
	The rank of the Jacobian matrix $Dg^{\lambda}$ is equal to $n-k$ for the affine subspace $V^{\lambda} = \mathbb{R}({^t (Q^{\lambda})^{-1}}p_{1})+ \cdots + \mathbb{R}({^t (Q^{\lambda})^{-1}}p_{k}) + {^t (Q^{\lambda})^{-1}}a$.

\begin{proof}
	As we have done in the beginning of this section,
	we choose a primitive basis $q_{k+1}, \ldots, q_{n}$ of the orthogonal subspace to $V$ spanned by the vectors $p_{1}, \ldots, p_{k}$.
	We can write $\overline{D_{O}(V)}$ as 
	\begin{equation*}
		\overline{D_{O}(V)}
		= \{\xi^{O} \in \Delta_{O} \mid g^{O}_{j} (\xi^{O}) = 0,\; j=k+1,\ldots,n \},
	\end{equation*}
	where the functions $g^{O}_{k+1},\ldots, g^{O}_{n}$ are defined by 
	\begin{equation*}
		g^{O}_{j}(\xi^{O})
		= 
		\prod_{i \in \mathcal{I}^{+}_{O,j}} (e \xi^{O}_{i})^{\langle e_{i},q_{j} \rangle} 
		- e^{\langle a,q_{j}\rangle} \prod_{i \in \mathcal{I}^{-}_{O,j}} (e \xi^{O}_{i})^{-\langle e_{i},q_{j} \rangle}
	\end{equation*}
	for $j=k+1,\ldots,n$. Here, the sets $\mathcal{I}^{+}_{O,j}$, $\mathcal{I}^{-}_{O,j}$ are defined by (\ref{eq: plus local}) and (\ref{eq: minus local}) respectively.
	We write the Jacobian matrix $Dg^{O} $ of $g^{O}= [g^{O}_{k+1}\; \cdots\; g^{O}_{n}]$ as 
	\begin{equation*}
		Dg^{O}
		= \left[
			\begin{matrix} \displaystyle
				\langle e_{j}, q_{i} \rangle e^{\langle a,q_{i} \rangle C^{O}_{ij}} \prod_{l \in \mathcal{I}_{O,i}}  e^{\pm \langle e_{l}, q_{i} \rangle} (\xi^{O}_{l})^{\pm \langle e_{l}, q_{i} \rangle - \delta_{lj}}
			\end{matrix}
		\right],
	\end{equation*}
	where $C^{O}_{ij}$ is defined by 
	\begin{equation*}
		C^{O}_{ij} =
		\begin{cases*}
			0 & if $j \in \mathcal{I}^{+}_{O,i}$, \\
			1 & if $j \not\in \mathcal{I}^{+}_{O,i}$.
		\end{cases*}
	\end{equation*}
	By the assumption, the rank of the Jacobian matrix $Dg^{O}$ of $g^{O}= [g^{O}_{k+1}\; \cdots\; g^{O}_{n}]$ is equal to $n-k$.

	We choose a primitive basis $q^{\lambda}_{k+1}, \ldots, q^{\lambda}_{n}$ of the orthogonal subspace to $V^{\lambda}$ spanned by the vectors ${^t (Q^{\lambda})^{-1}}p_{1}, \ldots, {^t (Q^{\lambda})^{-1}}p_{1}$ as 
	$q^{\lambda}_{k+1} = Q^{\lambda}q_{k+1}, \ldots, q^{\lambda}_{n} = Q^{\lambda}q_{n}$.
	Indeed, since $\langle p_{i},q_{j} \rangle = 0$ for any $i=1,\ldots,k$ and any $j=k+1,\ldots, n$, we have 
	\begin{align*}
		\left\langle {^t (Q^{\lambda})^{-1}}p_{i}, Q^{\lambda}q_{j} \right\rangle 
		&= {^t p_{i}}(Q^{\lambda})^{-1}Q^{\lambda}q_{j} \\
		&= \langle p_{i},q_{j} \rangle \\
		&= 0.
	\end{align*}
	We write $\overline{D_{\lambda}(V^{\lambda})}$ as 
	\begin{equation*}
		\overline{D_{\lambda}(V^{\lambda})}
		= \{\xi^{\lambda} \in \Delta_{\lambda} \mid g^{\lambda}_{j} (\xi^{\lambda}) = 0,\; j=k+1,\ldots, n \},
	\end{equation*}
	where the functions $g^{\lambda}_{k+1}, \ldots, g^{\lambda}_{n}$ are defined by 
	\begin{align*}
		g^{\lambda}_{j}(\xi^{\lambda})
		=& 
		\prod_{i\in \{\langle u^{\lambda}_{i}, Q^{\lambda}q_{j} \rangle \geq 0 \}} (e L^{\lambda}_{i})^{\langle u^{\lambda}_{i},Q^{\lambda}q_{j} \rangle} \\
		&- e^{\langle {^t (Q^{\lambda})^{-1}}a,Q^{\lambda}q_{j}\rangle}\prod_{i\in \{\langle u^{\lambda}_{i}, Q^{\lambda}q_{j} \rangle \leq 0 \}} (e L^{\lambda}_{i})^{-\langle u^{\lambda}_{i}, Q^{\lambda}q_{j} \rangle}
	\end{align*}
	for $j=k+1,\ldots, n$. 
	Recall that $Q^{\lambda} = \left[ v^{\lambda}_{1} \cdots v^{\lambda}_{n} \right]$, since 
	\begin{align*}
		\langle u^{\lambda}_{l}, Q^{\lambda} q_{i} \rangle 
		&= \langle {^t Q^{\lambda}}u^{\lambda}_{l}, q_{i} \rangle \\
		&= \left\langle \left[\begin{matrix} {^t v^{\lambda}_{1}} \\ \vdots \\ {^t v^{\lambda}_{n}} \end{matrix} \right] u^{\lambda}_{l}, q_{i} \right\rangle \\
		&= \left\langle \left[\begin{matrix} {^t v^{\lambda}_{1}} \\ \vdots \\ {^t v^{\lambda}_{n}} \end{matrix} \right] u^{\lambda}_{l}, q_{i} \right\rangle \\
		&= \langle e_{l}, q_{j} \rangle,
	\end{align*}
	and 
	\begin{align*}
		\langle {^t (Q^{\lambda})^{-1}}a,Q^{\lambda}q_{j}\rangle
		&=
		\langle {^t Q^{\lambda}}{^t (Q^{\lambda})^{-1}}a,q_{j}\rangle \\
		&=
		\langle a,q_{j}\rangle,
	\end{align*}
	we obtain 
	\begin{equation*}
		g^{\lambda}_{j}(\xi^{\lambda})
		= 
		\prod_{i\in \mathcal{I}^{+}_{O,j}} (e L^{\lambda}_{i})^{\langle e_{i},q_{j} \rangle} 
		- e^{\langle a,q_{j}\rangle}\prod_{i \in \mathcal{I}^{-}_{O,j}} (e L^{\lambda}_{i})^{-\langle e_{i}, Q^{\lambda}q_{j} \rangle},
	\end{equation*}
	for $j=k+1,\ldots,n$.
	We can calculate the Jacobian matrix $Dg^{\lambda}$ of $g^{\lambda}= [g^{\lambda}_{k+1}\; \cdots\; g^{\lambda}_{n}]$ as 
	\begin{align*}
		Dg^{\lambda} 
		&= \left[
			\begin{matrix} \displaystyle
				\sum_{m=1}^{n}\prod_{l \in \mathcal{I}^{\pm}_{O,i}} \langle e_{m}, q_{i} \rangle e^{\langle a,q_{i} \rangle C^{\lambda}_{ij}} \langle u^{\lambda}_{m},e_{j} \rangle e^{\pm \langle e_{l}, q_{i} \rangle} (L^{\lambda}_{l})^{\pm \langle e_{i}, q_{i} \rangle - \delta_{lm}}
			\end{matrix}
		\right] \\
		&= \left[
			\begin{matrix} \displaystyle
				\langle e_{j}, q_{i} \rangle e^{\langle a,q_{i} \rangle C^{O}_{ij}} \prod_{l \in \mathcal{I}^{\pm}_{O,i}} e^{\pm \langle e_{l}, q_{i} \rangle} (L^{\lambda}_{l})^{\pm \langle e_{l}, q_{i} \rangle - \delta_{lj}}
			\end{matrix}
		\right]
		\left[
			\begin{matrix} \displaystyle
				{^t u^{\lambda}_{1}} \\
				\vdots \\
				{^t u^{\lambda}_{n}}
			\end{matrix}
		\right].
	\end{align*}
	Thus, we obtain the expression of $Dg^{\lambda}$ as 
	\begin{equation*}
		Dg^{\lambda}
		= \left[
			\begin{matrix} \displaystyle
				\langle e_{j}, q_{i} \rangle e^{\langle a,q_{i} \rangle C^{O}_{ij}} \prod_{l \in \mathcal{I}^{\pm}_{O,i}} e^{\pm \langle e_{l}, q_{i} \rangle} (L^{\lambda}_{l})^{\pm \langle e_{l}, q_{i} \rangle - \delta_{lj}}
			\end{matrix}
		\right]
		\left[
			\begin{matrix} \displaystyle
				{^t u^{\lambda}_{1}} \\
				\vdots \\
				{^t u^{\lambda}_{n}}
			\end{matrix}
		\right].
	\end{equation*}
	The matrix $[u^{\lambda}_{1} \cdots u^{\lambda}_{n} ]$ does not influence the rank of $Dg^{\lambda}$.
	Lemma \ref{lemma: boundary of standard local to boundary of local} shows that $\xi^{O}_{i} = 0$ is equivalent to $L^{\lambda}_{i} = 0$ for $i =1,\ldots,n$.
	Comparing the matrix $Dg^{O}$ with $Dg^{\lambda}$, we see that 
	the rank of $Dg^{\lambda}$ is equal to the rank of $Dg^{O}$.
\end{proof}
\end{prop}

If $n-k = 1$, then it is sufficient to consider the linear part of a given $k$-dimensional affine subspace in $\mathbb{R}^{n}$ from Proposition \ref{prop: affine not affect when codim1} and Theorem \ref{thm: lyz correspondence}.
For a given Delzant polytope $\Delta$, define $\mathcal{V}_{\lambda}$ as the set of the $(n-1)$-dimensional subspaces $V$ in $\mathbb{R}^{n}$ such that the rank of the Jacobian matrix $Dg^{\lambda}$ is equal to one and define $\mathcal{V}$ as the set of the $(n-1)$-dimensional subspaces in $\mathbb{R}^{n}$ such that $\overline{D(V)}$ is a submanifold with corners in $\Delta$.
Then we obtain the set $\mathcal{V}$ by
\begin{equation}
	\mathcal{V} = \bigcap_{\lambda \in \Lambda} \mathcal{V}_{\lambda}.
\end{equation}
From Proposition \ref{prop: local models to global}, we see that the data of the standard local model $\mathcal{V}_{O}$ is sufficient to determine the set $\mathcal{V}$.

We give two examples for standard local models.
\begin{exam}
	\label{eg: standard local model n=2, k=1}
	For $n = 2$ and $k = 1$, we obtain 
	\begin{equation*}
		\mathcal{V}_{O}
		= \{ \langle p \rangle \mid p \in A_{++} \cup A_{+-} \},
	\end{equation*}
	where $\langle p \rangle$ is the subspace spanned by a vector $p$, and $A_{++}$ and $A_{+-}$ are defined by 
	\begin{align*}
		A_{++} &= 
		\left\{ 
			\begin{bmatrix}
				1 \\
				b_{2}
			\end{bmatrix} \in \mathbb{Z}^{2}
			\;\middle|\;
			b_{2} \geq 0
		\right\}
		\cup 
		\left\{ 
			\begin{bmatrix}
				b_{1} \\
				1
			\end{bmatrix} \in \mathbb{Z}^{2}
			\;\middle|\;
			b_{1} \geq 0
		\right\}, \\
		A_{+-} &=
		\left\{ 
			\begin{bmatrix}
				c_{1} \\
				c_{2}
			\end{bmatrix} \in \mathbb{Z}^{2}
			\;\middle|\;
			c_{1} \geq 0, c_{2} < 0
		\right\}.
	\end{align*}
\end{exam}

\begin{exam}
	For $n = 3$ and $k = 2$, we obtain
	\begin{equation*}
		\mathcal{V}_{O}
		= \{ \langle p \rangle \mid p \in B_{+++} \cup B_{-++} \cup B_{+-+} \cup B_{++-} \},
	\end{equation*}
	where $\langle p \rangle$ is the subspace spanned by a vector $p$, and $B_{+++}$, $B_{-++}$, $B_{+-+}$, and $B_{++-}$ are defined by
	\begin{align*}
		B_{+++} &=
		\left\{ 
			\begin{bmatrix}
				b_{1} \\ b_{2} \\ b_{3}
			\end{bmatrix} \in \mathbb{Z}^{3}
			\;\middle|\; b_{1} > 0, b_{2} > 0, b_{3} > 0 \right\}, \\
		B_{-++} &=
		\left\{
			\begin{bmatrix}
				-1 \\ b_{2} \\ b_{3}
			\end{bmatrix} \in \mathbb{Z}^{3}
			\;\middle|\; b_{2} \geq 0, b_{3} \geq 0 \right\}, \\
		B_{+-+} &=
		\left\{
			\begin{bmatrix}
				b_{1} \\ -1 \\ b_{3}
			\end{bmatrix} \in \mathbb{Z}^{3}
			\;\middle|\; b_{1} \geq 0, b_{3} \geq 0 \right\}, \\
		B_{++-} &=
		\left\{
			\begin{bmatrix}
				b_{1} \\ b_{2} \\ -1
			\end{bmatrix} \in \mathbb{Z}^{3}
			\;\middle|\; b_{1} \geq 0, b_{2} \geq 0 \right\}.
	\end{align*}
\end{exam}

In the rest of this section, 
we demonstrate the examples of submanifolds with corners in a Delzant polytope of $\mathbb{C}P^{2}$, which correspond to the examples of torus-equivariantly embedded toric manifolds in $\mathbb{C}P^{2}$ (see Remark \ref{remark: toric submanifold}).

From Lemma \ref{lemma: direction vectors vs normal vectors}, we see that ${^t (Q^{\lambda})^{-1}} = [\begin{matrix} u^{\lambda}_{1} & u^{\lambda}_{2} \end{matrix}]$ for $\lambda \in \Lambda = \{\lambda,\mu,\sigma \}$.
We define 
the points $\lambda$, $\mu$, $\sigma$ in $(\mathfrak{t}^{2})^{\ast} \cong \mathbb{R}^{2}$ by 
	\begin{equation*}
		\lambda = (0,0), ~
		\mu = (2,0), ~ 
		\sigma = (0,2).
	\end{equation*}
Let $\Delta$ be a polytope defined by the convex hull of the points $\lambda$, $\mu$, $\sigma$.
Then, we say that 
\begin{equation*}
	{^t (Q^{\lambda})^{-1}}
	= \left[
		\begin{matrix}
			1 & 0 \\
			0 & 1
		\end{matrix}
	\right], \;
	{^t (Q^{\mu})^{-1}}
	= \left[
		\begin{matrix}
			-1 & 1 \\
			-1 & 0
		\end{matrix}
	\right], \;
	{^t (Q^{\sigma})^{-1}}
	= \left[
		\begin{matrix}
			0 & -1 \\
			1 & -1
		\end{matrix}
	\right].
\end{equation*}
From Example \ref{eg: standard local model n=2, k=1}, we have obtained the standard model for $n = 2$ and $k = 1$.

Since $(Q^{\lambda})^{-1}$ is the identity matrix, we obtain the set $\mathcal{V}_{\lambda}$ as 
\begin{equation*}
	\mathcal{V}_{\lambda}
    =
    \{
        \langle p \rangle \mid p \in A^{\lambda}_{++}\cup A^{\lambda}_{+-}
    \},
\end{equation*}
where $\langle p \rangle$ is the subspace spanned by a vector $p$, and $A^{\lambda}_{++} = A_{++}$ and $A^{\lambda}_{+-} = A_{+-}$.
Since
\begin{align*}
	{^t(Q^{\mu})^{-1}} 
	\left[ \begin{matrix}
		1 \\ b_{2}	
	\end{matrix} \right]
	&= 
	\left[
		\begin{matrix}
			-1 + b_{2} \\
			-1
		\end{matrix}
	\right], \\
	{^t(Q^{\mu})^{-1}} 
	\left[ \begin{matrix}
		b_{1} \\ 1	
	\end{matrix} \right]
	&= 
	\left[
		\begin{matrix}
			-b_{1} + 1 \\
			-b_{1}
		\end{matrix}
	\right], \\
	{^t(Q^{\mu})^{-1}} 
	\left[ \begin{matrix}
		c_{1} \\ c_{2}	
	\end{matrix} \right]
	&= 
	\left[
		\begin{matrix}
			-c_{1} + c_{2} \\
			-c_{1}
		\end{matrix}
	\right],
\end{align*}
we obtain the set $\mathcal{V}_{\mu}$ as 
\begin{equation*}
	\mathcal{V}_{\mu}
	=
	\{
        \langle p \rangle \mid p \in A^{\mu}_{++}\cup A^{\mu}_{+-}
    \},
\end{equation*}
where $\langle p \rangle$ is the subspace spanned by a vector $p$, and $A^{\mu}_{++}$ and $A^{\mu}_{+-}$ are 
\begin{align*}
	A^{\mu}_{++} 
	&= 
	\left\{ 
		\left[
			\begin{matrix}
				b_{2} \\ -1 
			\end{matrix}
		\right] \in \mathbb{Z}^{2}
		\; \middle| \;
		b_{2} \geq -1
	\right\}
	\cup 
	\left\{
		\left[
			\begin{matrix}
				-b_{1} + 1 \\
				-b_{1}
			\end{matrix}
		\right] \in \mathbb{Z}^{2}
		\; \middle| \; 
		b_{1} \geq 0
	\right\}, \\
	A^{\mu}_{+-} 
	&=
	\left\{
		\left[
			\begin{matrix}
				-c_{1} + c_{2} \\
				-c_{1}
			\end{matrix}
		\right] \in \mathbb{Z}^{2}
		\; \middle| \;
			c_{1} \geq 0, \; c_{2} < 0
	\right\}.
\end{align*}
Similarly, since
\begin{align*}
	{^t(Q^{\sigma})^{-1}} 
	\left[ \begin{matrix}
		1 \\ b_{2}	
	\end{matrix} \right]
	&= 
	\left[
		\begin{matrix}
			- b_{2} \\
			1 - b_{2}
		\end{matrix}
	\right], \\
	{^t(Q^{\sigma})^{-1}} 
	\left[ \begin{matrix}
		b_{1} \\ 1	
	\end{matrix} \right]
	&= 
	\left[
		\begin{matrix}
			- 1 \\
			b_{1} - 1
		\end{matrix}
	\right], \\
	{^t(Q^{\sigma})^{-1}} 
	\left[ \begin{matrix}
		c_{1} \\ c_{2}	
	\end{matrix} \right]
	&= 
	\left[
		\begin{matrix}
			- c_{2} \\
			c_{1} - c_{2}
		\end{matrix}
	\right],
\end{align*}
we obtain the set $\mathcal{V}_{\sigma}$ as
\begin{equation*}
	\mathcal{V}_{\sigma}
	= 
    \{ 
        \langle p \rangle \mid p \in A^{\sigma}_{++} \cup A^{\sigma}_{+-} 
    \},
\end{equation*}
where $\langle p \rangle$ is the subspace spanned by a vector $p$, and $A^{\sigma}_{++}$ and $A^{\sigma}_{+-}$ are 
\begin{align*}
	A^{\sigma}_{++}
	&=
	\left\{ 
		\left[
			\begin{matrix}
				-b_{2} \\
				1 - b_{2}
			\end{matrix}
		\right] \in \mathbb{Z}^{2}
		\; \middle| \;
		b_{2} \geq 0
	\right\}
	\cup 
	\left\{
		\left[
			\begin{matrix}
				-1 \\
				b_{1} - 1
			\end{matrix}
		\right] \in \mathbb{Z}^{2}
		\; \middle| \;
		b_{1} \geq 0
	\right\}, \\
	A^{\sigma}_{+-}
	&=
	\left\{
		\left[
			\begin{matrix}
				-c_{2} \\
				c_{1} - c_{2}
			\end{matrix}
		\right] \in \mathbb{Z}^{2}
		\; \middle| \;
		c_{1} \geq 0,\; b_{2} < 0
	\right\}.
\end{align*}
We conclude that the set $\mathcal{V} = \mathcal{V}_{\lambda} \cap \mathcal{V}_{\mu} \cap \mathcal{V}_{\sigma}$ is expressed as 
\begin{equation*}
	\mathcal{V}
	= 
	\left\{
        \langle p \rangle 
        \; \middle| \; 
        p \in 
        \left\{
		\left[
			\begin{matrix}
				1 \\ 0
			\end{matrix}
		\right],
		\left[
			\begin{matrix}
				0 \\ 1
			\end{matrix}
		\right],
		\left[
			\begin{matrix}
				1 \\ 1
			\end{matrix}
		\right],
		\left[
			\begin{matrix}
				1 \\ -1
			\end{matrix}
		\right],
		\left[
			\begin{matrix}
				2 \\ 1
			\end{matrix}
		\right],
		\left[
			\begin{matrix}
				1 \\ 2
			\end{matrix}
		\right]
        \right\}
	\right\},
\end{equation*}
which coincides with the classification of torus-equivariantly embedded toric manifolds given in \cite[Remark 5.11]{yamaguchi2023toric}.

\section{Remarks on Intersection Theory}
\label{subsec: intersection theory}

We give some remarks on intersection theory for torus-equivariantly embedded toric manifolds $\overline{C(V)}$ and submanifolds $\overline{D(V)}$ with corners based on examples.

From our construction of $\mu: X \to (\mathfrak{t}^{n})^{\ast}$ and $\overline{D(V)}$, 
we see that 
\begin{equation*}
	\mu(\overline{C(V)})=\overline{D(V)}.
\end{equation*}
In \cite[Section 5]{yamaguchi2023toric}, examples for $\overline{C(V)}$ and corresponding $\overline{D(V)}$ were given.

In this section,
We concentrate on the case when $X = \mathbb{C}P^{2}$ and assume that $\overline{C(V)}$ is a complex submanifold in $\mathbb{C}P^{2}$.

\begin{figure}[h]
	\centering
	\begin{tabular}{cc}
	\begin{minipage}[t]{0.45\hsize}
		\centering
		\includegraphics[keepaspectratio,width=5cm]{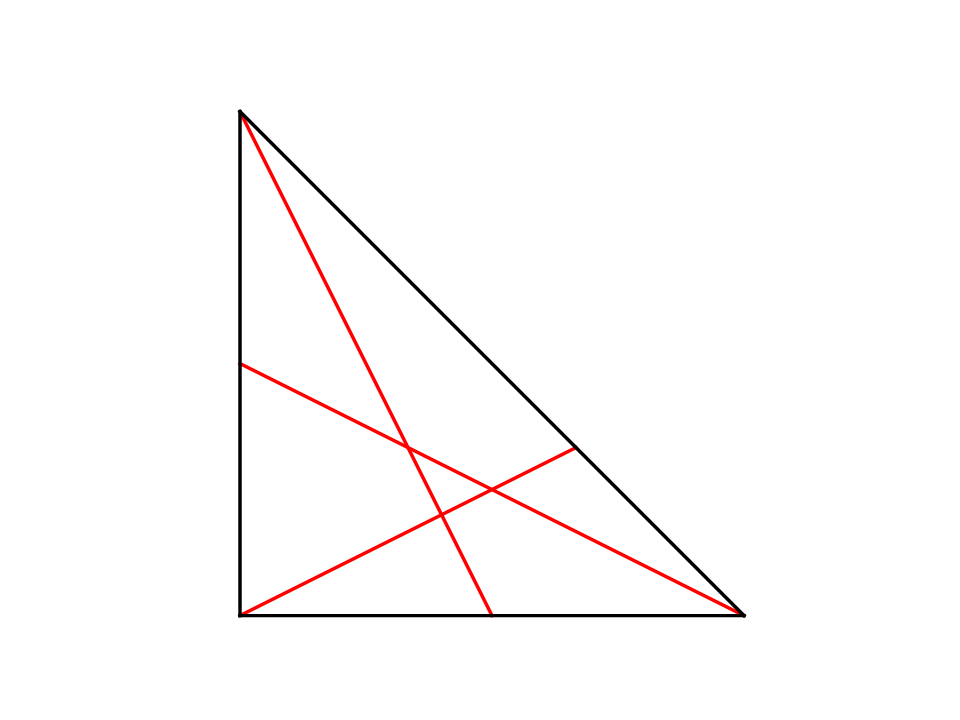}
		\caption{$\overline{D(V_{1})}$, $\overline{D(V_{2})}$, and $\overline{D(V_{3})}$}
		\label{fig: example100111a2CP2}
	\end{minipage} &
	\begin{minipage}[t]{0.45\hsize}
		\centering
		\includegraphics[keepaspectratio,width=5cm]{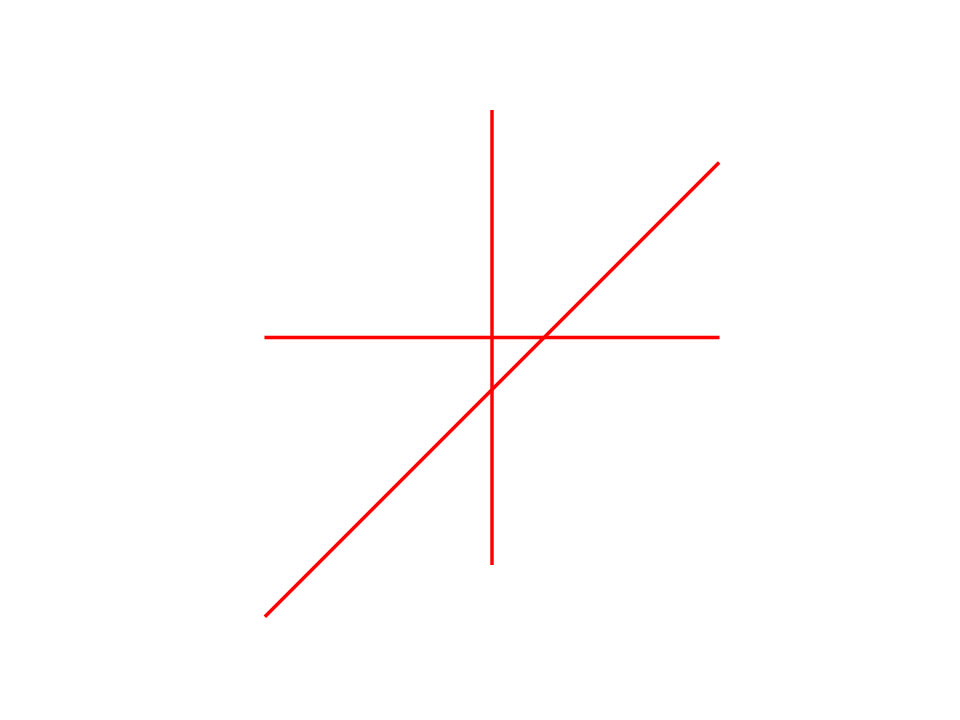}
		\caption{$V_{1}$, $V_{2}$, and $V_{3}$}
		\label{fig: example100111a2CP2 subspace}
	\end{minipage}
	\end{tabular}
\end{figure}

Let $V_{1}$, $V_{2}$, $V_{3}$ be affine subspaces defined by
\begin{equation*}
	V_{1} = \mathbb{R}\left[ \begin{matrix}
		1 \\ 0
	\end{matrix} \right], \;
	V_{2} = \mathbb{R} \left[ \begin{matrix}
		0 \\ 1
	\end{matrix} \right], \;
	V_{3} = \mathbb{R} \left[ \begin{matrix}
		1 \\ 1
	\end{matrix}\right] + \left[ \begin{matrix}
		\log 2 \\ 0
	\end{matrix}\right].
\end{equation*}
In Section \ref{subsec: classification of submanifolds in delzant}, we know that $\overline{C(V_{1})}$, $\overline{C(V_{2})}$, and $\overline{C(V_{3})}$ are complex submanifolds in $\mathbb{C}P^{2}$.
Figure \ref{fig: example100111a2CP2} describes the intersections of $\overline{D(V_{1})}$, $\overline{D(V_{2})}$, and $\overline{D(V_{3})}$ in a Delzant polytope of $\mathbb{C}P^{2}$. 
Figure \ref{fig: example100111a2CP2 subspace} describes the intersections of $V_{1}$, $V_{2}$, and $V_{3}$ in $\mathbb{R}^{2}$.

We notice that there is a one-to-one correspondence between the intersection points in the interior of the Delzant polytope of $\mathbb{C}P^{2}$ and the intersection points in $\mathbb{R}^{2}$.

\begin{figure}[h]
	\centering
	\begin{tabular}{cc}
	\begin{minipage}[t]{0.45\hsize}
		\centering
		\includegraphics[keepaspectratio,width=5cm]{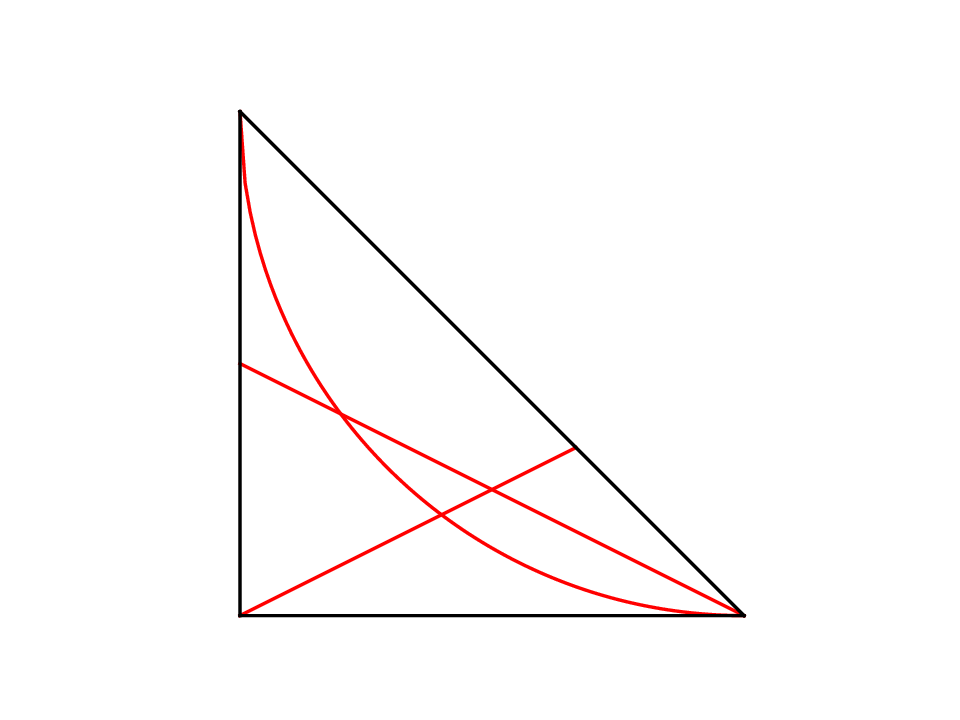}
		\caption{$\overline{D(V_{1})}$, $\overline{D(V_{3})}$, and $\overline{D(V_{4})}$}
		\label{fig: example11a21m1a110CP2}
	\end{minipage} &
	\begin{minipage}[t]{0.45\hsize}
		\centering
		\includegraphics[keepaspectratio,width=5cm]{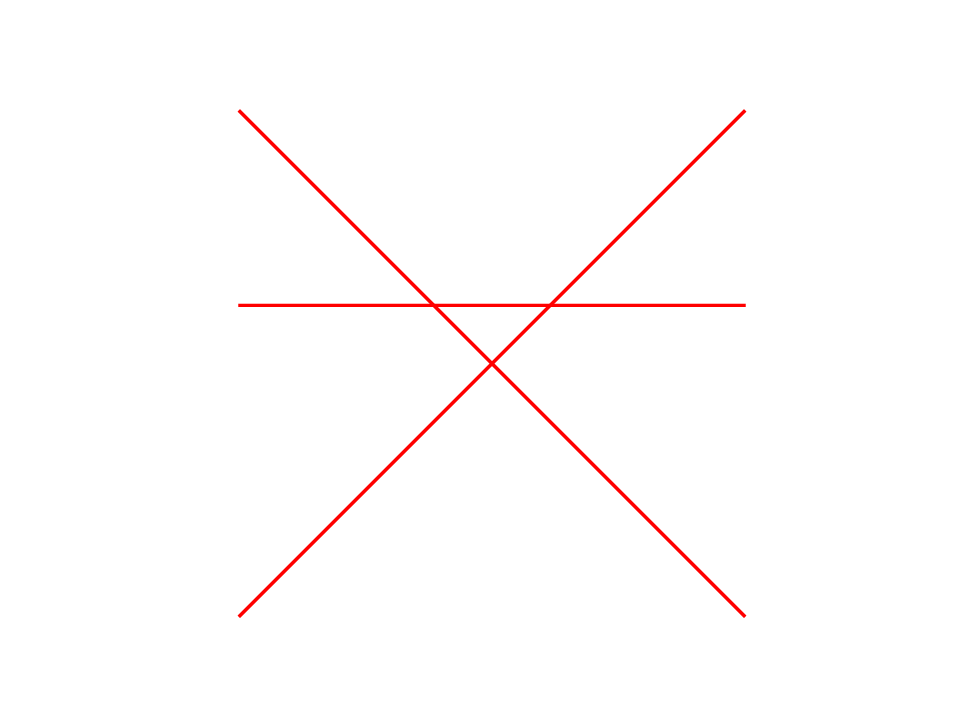}
		\caption{$V_{1}$, $V_{3}$, and $V_{4}$}
		\label{fig: example11a21m1a110CP2 subspace}
	\end{minipage}
	\end{tabular}
\end{figure}

Next, let $V_{4}$ be an affine subspace defined by
\begin{equation*}
	V_{4} = \mathbb{R} \left[ \begin{matrix}
		1 \\ -1
	\end{matrix}\right] + \left[ \begin{matrix}
		-\log 2 \\ 0
	\end{matrix}\right].
\end{equation*}
We know that $\overline{C(V_{4})}$ is also a complex submanifold in $\mathbb{C}P^{2}$.
Figure \ref{fig: example11a21m1a110CP2} describes the intersections of $\overline{D(V_{1})}$, $\overline{D(V_{3})}$, and $\overline{D(V_{4})}$ in a Delzant polytope of $\mathbb{C}P^{2}$. 
Figure \ref{fig: example11a21m1a110CP2 subspace} describes the intersections of $V_{1}$, $V_{3}$, and $V_{4}$ in $\mathbb{R}^{2}$.

We notice that $\overline{D(V_{1})}$ and $\overline{D(V_{4})}$ have an intersection point on the boundary of the Delzant polytope, whose corresponding intersection point cannot be found in $V_{1}$ and $V_{4}$.

\begin{figure}[h]
	\begin{center}
	\begin{tabular}{cc}
	\begin{minipage}[t]{0.45\hsize}
		\centering
		\includegraphics[keepaspectratio,width=5cm]{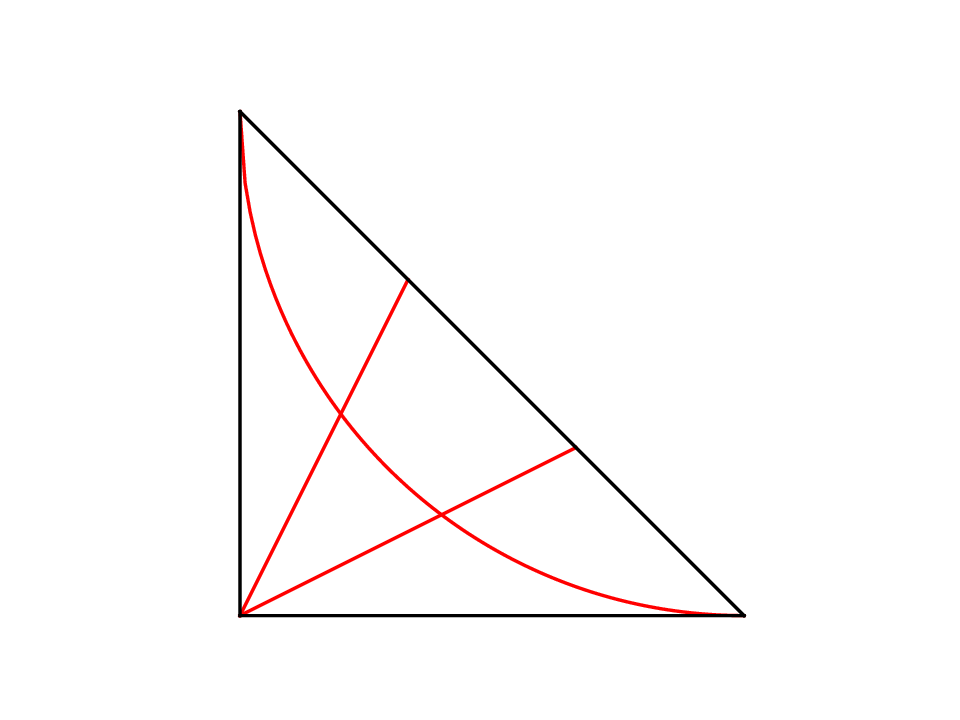}
		\caption{$\overline{D(V_{3})}$, $\overline{D(V_{4})}$, and $\overline{D(V_{5})}$}
		\label{fig: example11a111a21m1aCP2}
	\end{minipage} &
	\begin{minipage}[t]{0.45\hsize}
		\centering
		\includegraphics[keepaspectratio,width=5cm]{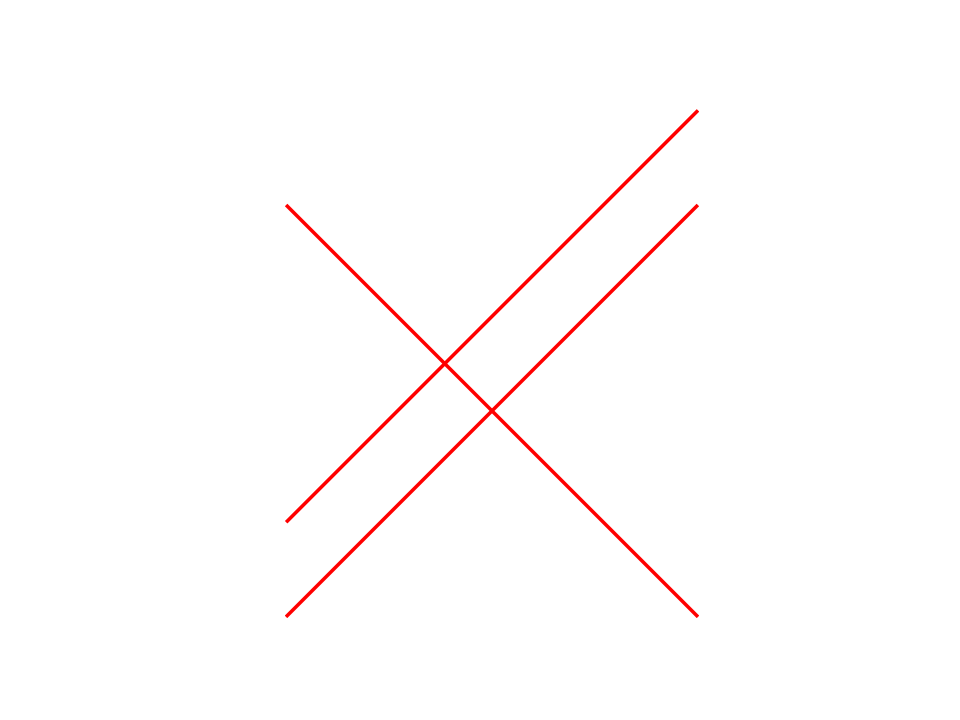}
		\caption{$V_{3}$, $V_{4}$, and $V_{5}$}
		\label{fig: example11a111a21m1aCP2 subspace}
	\end{minipage}
	\end{tabular}
    \end{center}   
\end{figure}

Furthermore, we give an example for this type of phenomena when one of three affine subspaces is parallel to one of the others.
Let $V_{5}$ be an affine subspace defined by 
\begin{equation*}
	V_{5} = \mathbb{R}\left[\begin{matrix}
		1 \\ 1
	\end{matrix} \right] + \left[ \begin{matrix}
		-\log 2 \\ 0
	\end{matrix} \right].
\end{equation*}
Since the $V_{5}$ is parallel to $V_{3}$, $\overline{C(V_{5})}$ is also a complex submanifold in $\mathbb{C}P^{2}$, as we know.
Figure \ref{fig: example11a111a21m1aCP2} describes the intersections of $\overline{D(V_{3})}$, $\overline{D(V_{4})}$, and $\overline{D(V_{5})}$ in a Delzant polytope of $\mathbb{C}P^{2}$. 
Figure \ref{fig: example11a111a21m1aCP2 subspace} describes the intersections of $V_{3}$, $V_{4}$, and $V_{5}$ in $\mathbb{R}^{2}$.

We notice that even if $V_{3}$ is parallel to $V_{5}$ in $\mathbb{R}^{2}$, $\overline{D(V_{3})}$, and $\overline{D(V_{5})}$ have an intersection point on the boundary of the Delzant polytope.

This observation might be useful when we consider Lagrangian intersection theory for the mirror Lagrangian submanifolds of torus-equivariantly embedded toric manifolds after making an explicit expression for the normal bundles of $\overline{D(V)}$ or the conormal bundles of $V$.
Moreover, we expect some relationship between intersection theory for torus-equivariantly embedded toric manifolds and intersection theory for tropical toric manifolds 
because some examples seem to be the same as in the context of tropical toric manifolds (see \cite{MR2428356} for example).

\section*{Acknowledgements}
The author is grateful to the advisor, Manabu Akaho for a lot of suggestions and supports. 
The author would like to thank Yuichi Kuno and Yasuhito Nakajima for helpful discussions, and would also like to thank the anonymous referee for helpful comments which improve the exposition.
This work is supported by JST, the establishment of university fellowships towards the creation of science technology innovation, Grant Number JPMJFS2139.

\bibliographystyle{alpha}
\bibliography{submanifolds_in_delzant}

\end{document}